\documentclass[a4paper,12pt]{article}
\usepackage{amsmath}
\usepackage{amsthm}
\usepackage{amssymb}
\usepackage{amscd}
\usepackage{hyperref}

\title{The Bochner identities for the K\"ahlerian gradients}
\author{Yasushi Homma 
\thanks{Department of Mathematical Sciences, Waseda University, 
  3-4-1 Ohkubo, Shinjuku-ku, Tokyo, 169-8555, JAPAN. \endgraf \textit{e-mail address}: homma@gm.math.waseda.ac.jp} }
 \date{}

\theoremstyle{plain}
\newtheorem{thm}{Theorem}[section]
\newtheorem{prop}[thm]{Proposition}
\newtheorem{cor}[thm]{Corollary}
\newtheorem{lem}[thm]{Lemma}
\theoremstyle{definition}
\newtheorem{defini}[thm]{Definition}
\newtheorem{exam}{Example}[section]

\numberwithin{equation}{section}

\theoremstyle{remark}
\newtheorem{rem}{Remark}[section]

\newcommand{\Ad}{\mathrm{Ad}}

\newcommand{\Hom}{\mathrm{Hom}}
\newcommand{\End}{\mathrm{End}}
\newcommand{\id}{\mathrm{id}}

\newcommand{\e}{\epsilon}

\newcommand{\bp}{\begin{pmatrix}}
\newcommand{\ep}{\end{pmatrix}}

\newcommand{\te}{ {}^t \! }
\newcommand{\gl}{\mathfrak{gl}(m,\mathbb{C})}

\allowdisplaybreaks[3]
\begin{document}
\maketitle
%%%%%%%%%%%%%%%%%%%%%%%%%%%%%%
%        abstract            %
%%%%%%%%%%%%%%%%%%%%%%%%%%%%%%
\begin{abstract}
We discuss algebraic properties for the symbols of geometric first order differential operators on almost Hermitian manifolds and K\"ahler manifolds. Through study on the universal enveloping algebra and higher Casimir elements, we know algebraic relations for the symbols like the Clifford algebra. From the relations, we have all the Bochner identities for the operators. As applications, we have vanishing theorems, the Bochner-Weitzenb\"ock formula, and eigenvalue estimates for the operators on K\"ahler manifolds.
2000MSC:
15A66, % Clifford algebras, spinors
17B35, % Universal enveloping (super)algebras
53A30, % Conformal differential geometry
53C55, % Hermitian and K\"ahlerian manifolds
58J05. % Elliptic equations on manifolds, general theory
\end{abstract}
\tableofcontents
%%%%%%%%%%%%%%%%%%%%%%%%%%%%%%%%
%%            0               %%
%%%%%%%%%%%%%%%%%%%%%%%%%%%%%%%%
\section{Introduction}\label{sec:0}
In Riemannian and spin geometry, geometric first order differential operators defined from the Riemannian metric $g$ (and spin structure) are important to research the structure of underlying manifold. The Dirac operator, the twistor operator, the exterior derivative, the interior derivative, and the conformal Killing operator are basic examples of geometric first order differential operators. These operators are known as conformally covariant first order differential operators, that is, under conformal deformation of the Riemannian metric, $g\mapsto g'=e^{2\sigma (x)}g$,  these operators are covariant. For example, the Dirac operator $D$ changes as
$$
D\mapsto D'=e^{(-\frac{n-1}{2}-1)\sigma}De^{\frac{n-1}{2}\sigma}.
$$
Here, the constant $\frac{n-1}{2}$ on this equation is called \textit{the conformal weight}, which depends on the highest weight of representation of $Spin(n)$. From the work of H. D. Fegan in \cite{Fe}, all the conformally covariant first order differential operators are realized as components of the covariant derivative $\nabla$ on associated vector bundles as follows. Let $M$ be a Riemannian manifold and $\mathbf{SO}(M)$ be the orthonormal frame bundle of $M$. The Levi-Civita connection induces the covariant derivative $\nabla$ on the associated bundle $\mathbf{S}_{\rho}:=\mathbf{SO}(M)\times_{\rho}V_{\rho}$, where $(\pi_{\rho},V_{\rho})$ is an irreducible unitary representation of $SO(n)$ with highest weight $\rho$. Then we have a first order differential operator defined by
 \begin{equation}
D^{\rho}_{\lambda}:\Gamma (\mathbf{S}_{\rho})\xrightarrow{\nabla}\Gamma (\mathbf{S}_{\rho}\otimes T^{\ast}(M))\xrightarrow{\Pi_{\lambda}}\Gamma (\mathbf{S}_{\lambda}),
\end{equation}
where $\Pi_{\lambda}$ is orthogonal projection onto irreducible component $\mathbf{S}_{\lambda}$ of $\mathbf{S}_{\rho}\otimes T^{\ast}(M)$. We call these operators \textit{the generalized gradients} or \textit{the Stein-Weiss gradients} (cf. \cite{SW}). We know that the gradients have conformal covariance, 
\begin{equation}
D^{\rho}_{\lambda}\mapsto D'{}^{\rho}_{\lambda}=e^{-(m(\rho,\lambda)+1)\sigma}D^{\rho}_{\lambda}e^{m(\rho,\lambda)\sigma}
\end{equation}
under conformal deformation. The constant $m(\rho,\lambda)$ is the conformal weight depending on highest weights $\rho$ and $\lambda$. Similaly, we have first order differential operators on spin manifolds. For example, we consider the spinor representation of $Spin(n)$ and have the Dirac operator and the twistor operator on the spinor bundle. 

In recent research for the gradients by T. Branson et al., we find out that the conformal weights are essential tool to investigate local properties of the gradients such as the Bochner identities, vanishing theorem, eigenvalue estimate, ellipticity, the Kato inequalities etc. (see \cite{Br2}, \cite{BH1}, \cite{BH2}, \cite{CGH}, \cite{H2}). In particular, some known results for the Dirac operator and the Laplace-Beltrami operator on differential forms can be obtained easily. The understanding of the gradients gives some new directions for global analysis and harmonic analysis as well as Riemannian and spin geometry. (see \cite{Br1}, \cite{Br3}, \cite{Br4}, \cite{Bu}, \cite{H1}, \cite{H3}).

In this paper, we change structure group $SO(n)$ or $Spin(n)$ into the unitary group $U(m)$. It means that we discuss the gradients on almost Hermitian manifolds. Let $M$ be an almost Hermitian manifold with real dimension $2m$ and $\mathbf{U}(M)$ be the principal bundle of unitary frames of $M$. Fix a connection on $\mathbf{U}(M)$, and we have a covariant derivative $\nabla$ on the vector bundle $\mathbf{S}_{\rho}=\mathbf{U}(M)\times_{\rho}V_{\rho}$, where $(\pi_{\rho},V_{\rho})$ is an irreducible unitary representation of $U(m)$. The covariant derivative splits as $\nabla=\nabla^{1,0}+\nabla^{0,1}$ with respect to almost complex structure. We decompose $\nabla^{1,0}$ and $\nabla^{0,1}$ further and have a set of first order differential operators as follows:
\begin{equation}
\begin{split}
D^{\rho}_{\lambda'}:\Gamma (\mathbf{S}_{\rho})\xrightarrow{\nabla^{1,0}} \Gamma (\mathbf{S}_{\rho}\otimes \Lambda^{1,0}(M))\xrightarrow{\Pi_{\lambda'}} \Gamma (\mathbf{S}_{\lambda'}), \\
D^{\rho}_{\lambda''}:\Gamma (\mathbf{S}_{\rho})\xrightarrow{\nabla^{0,1}} \Gamma (\mathbf{S}_{\rho}\otimes \Lambda^{0,1}(M))\xrightarrow{\Pi_{\lambda''}} \Gamma (\mathbf{S}_{\lambda''}).
\end{split}
\end{equation}
We call these operators \textit{the K\"ahlerian gradients}, which are the $U(m)$-invariant first order differential operators. Here, ``$U(m)$-invariant'' means that, if we have a diffeomorphism of $M$ preserving the Hermitian structure and the connection, then the operator is invariant. Moreover. we know that these operators also have conformal covariance under conformal deformation (see theorem \ref{thm:6-3}). 

Our aim is to give all the Bochner identities for the K\"ahlerian gradients (cf. \cite{Br2}, \cite{CGH}, \cite{H2} for the Riemannian case). Since the Bochner identities follow from local calculus of the operators, we discuss algebraic structure of the principal symbols of the operators like the Clifford algebra for the Dirac operator. We call the symbols \textit{the Clifford homomorphisms}, which are generalization of the Clifford multiplication on spinor spaces. The method to give the Bochner identities is as follows:
\begin{enumerate}
	\item We introduce the $(1,0)$- and $(0,1)$-conformal weights, which are K\"ahlerian analogue of the conformal weights in the Riemannian case. 
	\item By using the conformal weights, we relate the Clifford homomorphisms to the universal enveloping algebra $U(\gl)=U(\mathfrak{u}(m)\otimes \mathbb{C})$ (in theorem \ref{thm:3-5}). 
	\item Through investigation of $U(\gl)$, we have algebraic relations for the Clifford homomorphisms (in theorem \ref{thm:4-3} and corollary \ref{cor:4-5}). 	\item We lift the algebraic relations to the Bochner identities for the K\"ahlerian gradients (in theorem \ref{thm:7-4}).
	\end{enumerate}
We conclude that, in K\"ahler geometry, the $(0,1)$ and $(1,0)$-conformal weights are essential for local calculus of the K\"ahlerian gradients. Because many geometric operators in K\"ahler and spin K\"ahler geometry are realized as the K\"ahlerian gradients, our Bochner identities have a lot of applications to such geometry. Some vanishing theorems, the Bochner-Weitzenb\"ock formula, and estimates of eigenvalues are presented in this paper. We especially study about the case of the Dolbeault-Dirac operator and the Dirac operator. 

This paper is organized as follows. In section \ref{sec:1} and \ref{sec:2}, we prepare the representations, enveloping algebra, and (higher) Casimir elements of the unitary group or its Lie algebra. We introduce the $(1,0)$- and $(0,1)$-conformal weights, which appear naturally on spectral resolution of the Casimir elements. In section \ref{sec:3} and \ref{sec:4}, we introduce the notion of the Clifford homomorphisms as principal symbols of the K\"ahlerian gradients, and investigate algebraic structure of the Clifford homomorphisms. In section \ref{sec:5}, the spinors and the Clifford multiplications are given as an example. In section \ref{sec:6}, we define the K\"ahlerian gradients and show their conformal covariance. In section \ref{sec:7}, we gain the Bochner identities for the K\"ahlerian gradients. In the last section, some applications of the Bochner identities are presented. We know that many results in K\"ahler and spin K\"ahler geometry follows from our Bochner identities. 
%%%%%%%%%%%%%%%%%%%%%%%%%%%%%%%%
%%           1                %%
%%%%%%%%%%%%%%%%%%%%%%%%%%%%%%%%
\section{Representations of $U(m)$}\label{sec:1}
In this section, we give a short review to representation theory for unitary groups and their Lie algebras. Let $V$ be a $2m$-dimensional real vector space  with almost complex structure $J$ and inner product $g(\cdot,\cdot)$ such that $g(J u,Jv)=g(u,v)$ for any $u,v\in V$. The complexification $V\otimes \mathbb{C}$ of $V$ splits into the direct sum of the $(1,0)$-part $V^{1.0}$ and the $(0,1)$-part $V^{0,1}$ with respect to $J$. We know that $V$ is isomorphic to $V^{1,0}$ as a complex vector space and $V^{0,1}$ is complex conjugate to $V^{1,0}$, namely, $\overline{V^{1,0}}=V^{0,1}$. These $m$-dimensional complex vector spaces have natural Hermitian inner products as 
\begin{gather}
(u,v):=g(u,\bar{v})  \quad \textrm{for $u,v$ in $V^{1,0}$},  \nonumber\\
(\bar{u},\bar{v}):=g(\bar{u},v)  \quad \textrm{for $\bar{u},\bar{v}$ in $V^{0,1}$},  \nonumber
\end{gather}
where the inner product $g(\cdot,\cdot)$ on $V$ is extended complex linearly to the complex inner product on $V\otimes \mathbb{C}$. With these Hermitian inner products, we have isomorphisms $V^{1,0}\simeq (V^{0,1})^{\ast}$ and $V^{0,1}\simeq (V^{1,0})^{\ast}$. Let $e_1,Je_1,\cdots, e_m,Je_m$ be an orthonormal basis of $V$. Then Hermitian bases for $V^{1,0}$ and $V^{0,1}$ are 
\begin{equation}
\begin{cases}
  \e_k:=\frac{1}{\sqrt{2}}(e_k-\sqrt{-1}Je_k) \in V^{1,0}, \\
  \bar{\e}_k:=\frac{1}{\sqrt{2}}(e_k+\sqrt{-1}Je_k) \in V^{0,1} 
\end{cases}\label{eqn:1-1}
\end{equation}
for $k=1,\cdots m$. To simplify explanations, we set $V:=\mathbb{R}^{2m}$, $V^{1,0}:=\mathbb{C}^m$, $V^{0,1}:=\overline{\mathbb{C}^{m}}$, and choose the standard orthonormal basis of $\mathbb{R}^{2m}$. 

Let $U(m)$ be the unitary group on $\mathbb{C}^m$, and $\mathfrak{u}(m)$ be its Lie algebra. The complexification of $\mathfrak{u}(m)$ is the complex Lie algebra $\mathfrak{gl}(m,\mathbb{C})$ of $m\times m$ matrices with real structure $Z\to -\te\bar{Z}$. For $k,l=1,\cdots,m$, we define a matrix $e_{kl}$ to be $1$ in the $(k,l)$-th place and $0$ elsewhere. Then $\{e_{kl}\}_{1\le k,l\le m}$ constitute a basis of $\mathfrak{gl}(m,\mathbb{C})$ and satisfy 
\begin{equation}
[e_{ij},e_{kl}]=\delta_{jk}e_{il}-\delta_{li}e_{kj} \quad \textrm{for any $i,j,k,l$}. \label{eqn:1-2}
\end{equation}

We choose $\mathfrak{h}=\mathrm{span}_{\mathbb{R}}\{ \sqrt{-1}e_{kk}|1\le k \le m \}$ as a cartan subalgebra, which is a maximal abelian subalgebra of $\mathfrak{u}(m)$. On a finite-dimensional unitary representation $(\pi,V)$ of $\mathfrak{u}(m)$, we can decompose the representation space into simultaneous eigenspaces called weight spaces with respect to $\mathfrak{h}$. For each weight space, we have an $m$-dimensional vector $\lambda=(\lambda^1,\cdots,\lambda^m)$ such that $\lambda^k$ is the eigenvalue of $e_{kk}$. This $\lambda$ is called the weight of weight space, and each component $\lambda^k$ is an integer. Now, for an irreducible representation $(\pi,V)$, we order the weights lexicographically, and get the highest one $\rho$. The highest weight $\rho$ satisfies \textit{the dominant integral condition}, 
%%%%%%%% weight %%%%%%%
\begin{equation}
\rho=(\rho^1,\cdots,\rho^m)\in \mathbb{Z}^m \quad \textrm{and} \quad \rho^1\ge \cdots \ge \rho^{m-1}\ge \rho^m. \label{eqn:1-3}
\end{equation}
%%%%%%
Conversely, for a vector $\rho$ with dominant integral condition, we have an irreducible unitary representation of $\mathfrak{u}(m)$ with highest weight $\rho$. Moreover, there is one-to-one correspondence between the finite-dimensional representations of $U(m)$ and its infinitesimal ones of $\mathfrak{u}(m)$. So we denote by $(\pi_{\rho},V_{\rho})$ an irreducible unitary representation of $U(m)$ or $\mathfrak{u}(m)$ with highest weight $\rho$. 

Let us see some irreducible $U(m)$-modules. When writing weights, we denote a string of $p$ $k$'s for $k$ in $\mathbb{Z}$ by $k_p$. For example, $(1_p,0_{m-p})$ is the weight whose first $p$ components are $1$ and others are $0$. 
%%%%%%%%% example %%%%%%%%%%%%%%
\begin{exam}\label{ex:1-1}
The natural representation of $U(m)$ on $V^{1,0}=\mathbb{C}^m$ has the highest weight $(1,0_{m-1})$. The conjugate representation on $V^{0,1}=\overline{\mathbb{C}^m}$ or the contragredient representation $(V^{1,0})^{\ast}\simeq V^{0,1}$ has the highest weight $(0_{m-1},-1)$.
\end{exam}
%%%%%%%%% example %%%%%%%%%%%%%%
\begin{exam}\label{ex:1-2}
The $p$-th exterior tensor product representation on $\Lambda^p V^{1,0}=\Lambda^p \mathbb{C}^m$ has the highest weight $(1_p,0_{m-p})$. Its conjugate representation on $\Lambda^p V^{0,1}=\Lambda^p\overline{\mathbb{C}^m}$ has the highest weight $(0_{m-p},(-1)_p)$. The $p$-th symmetric tensor product representation on $S^p(\mathbb{C}^m)$ has the highest weight $(p,0_{m-1})$.
\end{exam}
%%%%%%%%% example %%%%%%%%%%%%%%
\begin{exam}\label{ex:1-3}
The contragredient representation $(\pi_{\te \rho},V_{\te \rho})$ of $(\pi_{\rho},V_{\rho})$ has the highest weight $\te\rho=(-\rho^m,-\rho^{m-1},\cdots, -\rho^1)$.
\end{exam}
%%%%%%%%%%%%%%%%%%%%%%%%%%%%%%%%%%%%%
%                2                  %
%%%%%%%%%%%%%%%%%%%%%%%%%%%%%%%%%%%%%
\section{Casimir elements for $U(m)$}\label{sec:2}
The unitary representations of $\mathfrak{u}(m)$ correspond to the complex representations of $\gl$ and of the universal enveloping algebra $U(\gl)$. So we denote the extension of the representation $(\pi_{\rho}, V_{\rho})$ to $\gl$ and $U(\gl)$ by the same notation $(\pi_{\rho},V_{\rho})$. 

The universal enveloping algebra $U(\mathfrak{gl}(m,\mathbb{C}))$ of $\mathfrak{gl}(m,\mathbb{C})$ is the quotient algebra of the tensor algebra $T(\mathfrak{gl}(m,\mathbb{C}))$  by the two-sided ideal generated by all $(X\otimes Y-Y\otimes X-[X,Y])$ for $X,Y$ in $\mathfrak{gl}(m,\mathbb{C})$. When we define $U^N(\gl)$ to be the image of $\sum_{r=0}^N\otimes^r \gl$, we know that the algebra $U(\mathfrak{gl}(m,\mathbb{C}))$ has a filtered algebra structure. If an element $z$ is in $U^q(\gl)$, but is not in $U^{q-1}(\gl)$, then the degree of $z$ is said to be $q$.

We shall discuss the center $\mathfrak{Z}$ and the (higher) Casimir elements of $U(\mathfrak{gl}(m,\mathbb{C}))$. For every non-negative integer $q$, we define an element $e_{kl}^q$ with degree $q$ by 
\begin{equation}
e_{kl}^q:=
\begin{cases}
\sum_{ 1\le i_1,i_2,\cdots i_{q-1}\le m } e_{ki_1} e_{i_1i_2} \cdots e_{i_{q-1}l},  & \textrm{for $q \ge 1$}, \\
\delta_{kl} & \textrm{for $q=0$}. 
\end{cases}\label{eqn:2-1}
\end{equation}
%%%%%%%%%   lemma  %%%%%%%%%%%%%%%%%%%%
\begin{lem}\label{lem:2-1}
The elements $\{e_{kl}^q|q\in \mathbb{Z}_{\ge 0}, k,l=1,\cdots,m\}$ satisfy 
\begin{gather}
[e_{ij},e_{kl}^q]=\delta_{jk}e_{il}^q-\delta_{il}e^q_{kj}, \label{eqn:2-2} \\
\sum_{1\le i\le m} e_{ki}^p e_{il}^q=e_{kl}^{p+q}. \label{eqn:2-3}
\end{gather}
\end{lem}
%%%%%%%%% proof %%%%%%%%
\begin{proof}
The second equation is clear. Let us prove the first equation. The equation \eqref{eqn:2-2} is trivial for $q=0, 1$. For the general case, we consider the adjoint representation of $\gl$ on $U(\gl)$ and use the equation \eqref{eqn:1-2}.
\end{proof}
%%%%
From \eqref{eqn:2-2}, the trace $c_q:=\sum_k e_{kk}^q$ is a Casimir element with degree $q$, that is, $c_q$ is in $\mathfrak{Z}\cap U^q(\gl)$. The following fact for $c_q$ is well-known (see \cite{Z}). 
%%%%%%%%%%%%%%%%%  proposition %%%%%%%%%%%%%%%%%%%
\begin{prop}[\cite{Z}]\label{prop:2-2}
Let $c_q$ be the Casimir element defined by $c_q:=\sum_{1\le k\le m} e_{kk}^q$ for $q$ in $\mathbb{Z}_{\ge 0}$. On irreducible $U(m)$-module $V_{\rho}$ with highest weight $\rho=(\rho^1,\cdots,\rho^m)$, the Casimir element $c_q$ is constant as follows:
	\begin{equation}
	\pi_{\rho}(c_q)=\sum_{i=1}^m w_{-i}^q\gamma_{-i}, \label{eqn:2-4}
	\end{equation}
	where 
	\begin{equation}
	w_{-i}:=\rho^i+(m-i), \quad \gamma_{-i}:=\prod_{j\neq i}\left( 1-\frac{1}{w_{-i}-w_{-j}}\right) \label{eqn:2-5}
	\end{equation}	 
	for $1\le i \le m$. 
\end{prop}
%%%%%
We call the above constant $w_{-i}$ \textit{the $(0,1)$-conformal weight associated to $\rho$}.

The Casimir elements with the first few degree have simple descriptions.
\begin{enumerate}
	\item The $0$-th Casimir element is $c_0=\sum_k \delta_{kk}=m$. So we have $\pi_{\rho}(c_0)=m$ for any $\rho$. 
	\item The first Casimir element is $c_1=\sum_k e_{kk}$. So we have $\pi_{\rho}(c_1)=\sum_{i} \rho^i$. 
	\item The second Casimir element is the usual Casimir element $c_2=\sum_{kl} e_{kl}e_{lk}$. Then we show that
	\begin{equation}
\pi_{\rho}(c_2)=\sum_{k,l} \pi_{\rho}(e_{kl})\pi_{\rho}(e_{lk})=\sum_i \rho^i(\rho^i+m-2i+1). \label{eqn:2-6}
	\end{equation}
\end{enumerate}

We shall introduce new Casimir elements. The Lie algebra $\gl$ has the involution\begin{equation}
\gl \ni Z\mapsto \tilde{Z}:=-\te Z\in \gl \nonumber 
\end{equation}
such that $\widetilde{[Z,W]}=[\tilde{Z}, \tilde{W}]$. This involution naturally extends to an algebraic automorphism of $U(\gl)$, 
\begin{multline}
 U^q(\gl)\ni \omega=Z_1\cdots Z_q \\ \mapsto \tilde{\omega}=(-1)^q \: \te Z_1\cdots \te Z_q \in U^q(\gl). \label{eqn:2-8}
\end{multline}
We remark that this involution preserves the center $\mathfrak{Z}$, namely, $\tilde{\mathfrak{Z}}=\mathfrak{Z}$. So $\tilde{c}_q$ is also a Casimir element. To investigate such Casimir elements, we consider $\tilde{e}_{kl}^q:=\widetilde{e_{kl}^q}$, 
\begin{equation}
\tilde{e}_{kl}^q=(-1)^q\sum_{1\le i_1,\cdots,i_{q-1}\le m}e_{i_1k}e_{i_2i_1}\cdots e_{li_{q-1}} \in U^q(\gl). \label{eqn:2-9}
\end{equation}
The involutions of \eqref{eqn:2-2} and \eqref{eqn:2-3} lead us to the next lemma. 
%%%%%%%%%%% lemma %%%%%%%%%%%%%%
\begin{lem}\label{lem:2-3}
The elements $\{\tilde{e}_{kl}^q|q\in\mathbb{Z}_{\ge 0}, k,l=1,\cdots,m\}$ satisfy 
\begin{gather}
[e_{ij},\tilde{e}_{kl}^q]=\delta_{jl}\tilde{e}_{ki}^q-\delta_{ik}\tilde{e}^q_{jl}, \label{eqn:2-10}\\
\sum_i \tilde{e}_{ki}^p \tilde{e}_{il}^q=\tilde{e}_{kl}^{p+q}. \label{eqn:2-11}
\end{gather}
\end{lem}
%%%%%%
We calculate the eigenvalue of $\tilde{c}_q=\sum \tilde{e}_{kk}^q$ on irreducible $U(m)$-module $V_{\rho}$. It is from the definition of $\tilde{c}_q$ that
\begin{equation}
\begin{split}
\pi_{\rho}(\tilde{c}_q)& =(-1)^q\sum_{1\le i_1,\cdots,i_{q-1},k\le m} \pi_{\rho}(e_{i_1k}) \pi_{\rho}(e_{i_2i_1}) \cdots \pi_{\rho}(e_{ki_{q-1}})\\ 
 &=\sum_{1\le i_1,\cdots,i_{q-1},k \le m} \pi_{\te \rho}(e_{ki_1}) \pi_{\te \rho}(e_{i_1i_2}) \cdots \pi_{\te \rho}(e_{i_{q-1}k}) \\
&=\pi_{\te \rho}(c_q),
\end{split}\nonumber
\end{equation}
where $\pi_{\te \rho}$ is the contragredient representation of $\pi_{\rho}$ whose highest weight is $\te \rho=(-\rho^m,\cdots,-\rho^1)$. 
So we have 
$$
\pi_{\rho}(\tilde{c}_q)=\pi_{\te \rho}(c_q)=\sum_i (\te w_{-i})^q \: \te \gamma_{-i}.
$$
Here, 
\begin{equation}
\te w_{-i}:=(\te \rho)^i+(m-i)=-\rho^{m-i+1}+(m-i),  \quad \te \gamma_{-i}:=\prod_{j\neq i}\left(1-\frac{1}{\te w_{-i}-\te w_{-j} }\right). \nonumber 
\end{equation}
When we set 
\begin{equation}
w_{+i}:=-\rho^i+i-1, \quad \gamma_{+i}:=\prod_{j\neq i}\left(1-\frac{1}{w_{+i}- w_{+j} }\right) \label{eqn:2-12}
\end{equation}
 for $i=1,\cdots m$, we show that $\sum_i ({}^t\! w_{-i})^q \: \te \gamma_{-i}=\sum_i w_{+i}^q\gamma_{+i}$. Hence, we conclude that 
%%%%%%%%%%%%%%%% proposition %%%%%%%%%%%%
\begin{prop}\label{prop:2-4}
The Casimir element $\tilde{c}_q$ is the following constant on irreducible $U(m)$-module $V_{\rho}$:
\begin{equation}
\pi_{\rho}(\tilde{c}_q)=\sum_i w_{+i}^q \gamma_{+i}, \label{eqn:2-13}
\end{equation}
where $w_{+i}$ and $\gamma_{+i}$ are given by \eqref{eqn:2-12}. 
\end{prop}
We call the constant $w_{+i}$ \textit{the $(1,0)$-conformal weight associated to $\rho$}. 

%%%%%%%%%%%%%%%%%%%%%%%%%%%%%%%%%%%%%%%%%%%%%%%%%%%%%%%%%%%%%%
%%%%%%%%%%%%%%%%%%%%%%%%%%%%%%%%
%%            3               %%
%%%%%%%%%%%%%%%%%%%%%%%%%%%%%%%%
\section{Clifford homomorphisms for $U(m)$}\label{sec:3}
In this section, we define a generalization of the Clifford multiplication corresponding to the symbols of the K\"ahlerian gradients. First, we consider the irreducible unitary representation $(\pi_{\rho},V_{\rho})$ and the tensor product representation $(\pi_{\rho}\otimes \pi_{\mu_1},V_{\rho}\otimes \mathbb{C}^m)$. Here, $\{\mu_i\}_{i=1}^m$ is the standard basis of $\mathbb{Z}^m$,
\begin{equation}
\mu_i:=(0_{i-1},1,0_{m-i})=(\underbrace{0,\cdots, 0}_{i-1},1,\underbrace{0,\cdots, 0}_{m-i})\in \mathbb{Z}^m.  \label{eqn:3-1}
\end{equation}
By using the tensor product decomposition rules (so-called the Littlewood-Richardson rules) in \cite{Ko}, \cite{Z}, we show that highest weights of the irreducible components of $V_{\rho}\otimes \mathbb{C}^m$ are 
$$
\{\rho+\mu_i| \textrm{ $\rho+\mu_i$ is dominant integral, } 1\le i\le m \}.  
$$
When $\rho+\mu_i$ does  not satisfy the dominant integral condition, we set $V_{\rho+\mu_i}:=\{0\}$. Thus, we have the decomposition 
\begin{equation}
V_{\rho}\otimes \mathbb{C}^m=\sum_{1\le i\le m}V_{\rho+\mu_i}. \label{eqn:3-3}
\end{equation}
Next, we consider the tensor product representation $(\pi_{\rho}\otimes \pi_{-\mu_m},V_{\rho}\otimes \overline{\mathbb{C}^m})$. Then highest weights of irreducible components of $V_{\rho}\otimes \overline{\mathbb{C}^m}$ are
$$
\{\rho-\mu_i | \textrm{ $\rho-\mu_i$ is dominant integral, } 1\le i\le m\}.
$$
So we set $V_{\rho-\mu_i}:=\{0\}$ for $\rho-\mu_i$ without dominant integral condition and have the irreducible decomposition 
\begin{equation}
V_{\rho}\otimes \overline{\mathbb{C}^m}=\sum_{1\le i\le m}V_{\rho-\mu_i}. \label{eqn:3-5}
\end{equation}
The following is the definition of the Clifford homomorphism introduced by the author in \cite{H1}, \cite{H2}, which is a generalization of the Clifford multiplications on spinor space. 
%%%%%%%%%%%%%% definition %%%%%%%%%%%%%%%%%
\begin{defini}\label{def:3-1}
\begin{enumerate}
	\item Let $\Pi_{+i}$ be the orthogonal projection from $V_{\rho}\otimes \mathbb{C}^m$ onto $V_{\rho+\mu_i}$. For $u$ in $\mathbb{C}^m$, we define the linear mapping $p_{+i}(u)$  from $V_{\rho}$ to $V_{\rho+\mu_i}$ by 
\begin{equation}
p_{+i}(u):V_{\rho}\ni \phi \mapsto \Pi_{+i}(\phi\otimes u) \in V_{\rho+\mu_i}. \label{eqn:3-6}
\end{equation}
We denote by $p_{+i}(u)^{\ast}$ the adjoint of $p_{+i}(u)$ such that $(p_{+i}(u)\phi,\psi)=(\phi,p_{+i}(u)^{\ast}\psi)$ for $\phi$ in $V_{\rho}$ and $\psi$ in $V_{\rho+\mu_i}$.
\item Let $\Pi_{-i}$ be the orthogonal projection from $V_{\rho}\otimes \overline{\mathbb{C}^m}$ onto $V_{\rho-\mu_i}$. For $\bar{u}$ in $\overline{\mathbb{C}^m}$, we define the linear mapping 
\begin{equation}
p_{-i}(\bar{u}):V_{\rho}\ni \phi \mapsto \Pi_{-i}(\phi\otimes \bar{u}) \in V_{\rho-\mu_i}, \label{eqn:3-7}
\end{equation}
and denote its adjoint by $p_{-i}(\bar{u})^{\ast}$. 
\end{enumerate}
We call $p_{\pm i}$ and $p_{\pm i}^{\ast}$ \textit{the Clifford homomorphisms}.  We sometimes denote $p_{\pm i}$ by $p^{\rho}_{\rho\pm \mu_i}$ to clarify considered representation spaces.
\end{defini}
%%%%%%%%%%%% remark %%%%%%%%%%%%
\begin{rem}\label{rem:3-2}
Because multiplicity of each irreducible component is one, the orthogonal projection $\Pi_{\pm i}$ is well-defined. But, when we realize the representation spaces concretely, the Clifford homomorphisms have an ambiguity of complex number $a$ with $|a|=1$. For example, we consider the decomposition $\mathbb{C}^m\otimes \overline{\mathbb{C}^m}=V_{(0_m)}\oplus V_{(1,0_{m-2},-1)}$. For $\e_k\otimes \bar{\e}_l$ in $\mathbb{C}^m\otimes \overline{\mathbb{C}^m}$, we have
$$
\e_k\otimes \bar{\e}_l=\frac{\delta_{kl}}{n}\sum \e_i\otimes \bar{\e}_i+(-\frac{\delta_{kl}}{n}\sum \e_i\otimes \bar{\e}_i+\e_k\otimes \bar{\e}_l)\in V_{(0_m)}\oplus V_{(1,0_{m-2},-1)}.
$$
In this case, the Clifford homomorphism $p_{-1}(\bar{\e}_l)$ is just the interior product $i(\bar{\e}_l)$,
$$
i(\bar{\e}_l):\mathbb{C}^m \ni \e_k\mapsto \delta_{kl}\in \mathbb{C},
$$
where we use the identification
$$
V_{(0_m)}\ni \frac{1}{n}\sum \bar{\e}_i\otimes \e_i \mapsto 	1 \in \mathbb{C}.
$$
If we choose another identification preserving inner product, $\frac{1}{n}\sum \bar{\e}_i\otimes \e_i \mapsto e^{-i\theta}$, then the Clifford homomorphism is realized as $e^{-i\theta}i(\bar{\e}_k)$.
\end{rem}
%%%%

We shall provide some basic properties of the Clifford homomorphisms. First, we have
%%%%%%%%%%%%%% proposition %%%%%%%%%%%%%%%
\begin{prop}\label{prop:3-2}
The Clifford homomorphisms satisfy that, for $g$ in $U(m)$ and $u$ in $\mathbb{C}^m$, 
\begin{equation}
\begin{split}
p_{+i}(gu)&=\pi_{\rho+\mu_i}(g)p_{+i}(u)\pi_{\rho}(g^{-1}) :V_{\rho}\to V_{\rho+\mu_i} , \\
p_{+i}(gu)^{\ast}&=\pi_{\rho}(g)p_{+i}(u)^{\ast}\pi_{\rho+\mu_i}(g^{-1}) :V_{\rho+\mu_i}\to V_{\rho}, \\
p_{-i}(\overline{gu})&=\pi_{\rho-\mu_i}(g)p_{-i}(\bar{u})\pi_{\rho}(g^{-1}) :V_{\rho}\to V_{\rho-\mu_i},  \\
 p_{-i}(\overline{gu})^{\ast}&=\pi_{\rho}(g)p_{-i}(\bar{u})^{\ast}\pi_{\rho-\mu_i}(g^{-1}) :V_{\rho-\mu_i}\to V_{\rho}. 
\end{split} \label{eqn:3-8}
\end{equation}
Here, $gu$ means the natural representation $\pi_{\mu_1}(g)u$.
\end{prop}
%%%%%%%% proof %%%%%%%%%
\begin{proof}
Considering the action of $g$ to $\phi\otimes u$ in $V_{\rho}\otimes \mathbb{C}^m$, we have
$$
\sum_i \pi_{\rho+\mu_i}(g)p_{+i}(u)\phi =\pi_{\rho}\otimes\pi_{\mu_1}(g)(\phi\otimes u)=\pi_{\rho}(g)\phi\otimes gu 
         =\sum_i p_{+i}(gu)\pi_{\rho}(g)\phi.
$$
Then we have $\pi_{\rho+\mu_i}(g)p_{+i}(u)=p_{+i}(gu)\pi_{\rho}(g)$ for any $i$. \end{proof}
%%%%

Next, we show some algebraic relations of the Clifford homomorphisms. 
%%%%%%%%%%%%%% lemma %%%%%%%%%%%%%%
\begin{lem}\label{lem:3-3}
Let $\{ \e_k\}_k$ be a unitary basis of $\mathbb{C}^m$ and $\{\bar{\e}_k\}_k$ be a dual basis of $\overline{\mathbb{C}^m}$ given in \eqref{eqn:1-1}. Then
\begin{equation}
\sum_{1\le i\le m} p_{+i}(\e_k)^{\ast}p_{+i}(\e_l)=\delta_{kl}, \quad \sum_{1\le i\le m} p_{-i}(\bar{\e}_k)^{\ast}p_{-i}(\bar{\e}_l)=\delta_{kl} \label{eqn:3-10}
\end{equation}
for any $k,l$. 
\end{lem}
%%%%%%%%%%% proof %%%%%%%%%%%
\begin{proof}
From the definition of the Clifford homomorphisms, we show that
$$
\delta_{kl}(\phi,\psi)=(\phi\otimes \e_l,\psi\otimes \e_k)=\sum_i(p_{+i}(\e_l)\phi,p_{+i}(\e_k)\psi)=\sum_i(p_{+i}(\e_k)^{\ast}p_{+i}(\e_l)\phi,\psi)
$$
for any $\phi$, $\psi$ in $V_{\rho}$. Then we have proved the lemma. 
\end{proof}
%%%%%%%
To get more relations, we make use of the conformal weights $w_{\pm i}$ in \eqref{eqn:2-5} and \eqref{eqn:2-12}. Let us consider the tensor product representation $(\pi_{\rho}\otimes \pi_{-\mu_m},V_{\rho}\otimes \overline{\mathbb{C}^m})$ and define the operator $\widehat{C}_-$ on $V_{\rho}\otimes \overline{\mathbb{C}^m}$ by
\begin{equation}
\begin{split}
\widehat{C}_-: &=\pi_{\rho}\otimes \pi_{-\mu_m}(c_2)-\pi_{\rho}(c_2)\otimes \id -\id \otimes \pi_{-\mu_m}(c_2) \\
      &=2\sum_{k,l} \pi_{\rho}(e_{kl})\otimes \pi_{-\mu_m}(e_{lk}),
\end{split}\label{eqn:3-11}
\end{equation}
where $c_2$ is the Casimir element with degree $2$. This operator $\widehat{C}_-$ acts as a constant on each irreducible component $V_{\rho-\mu_i}$. In fact, it follows from \eqref{eqn:2-6} that 
\begin{equation}
\widehat{C}_-=-2w_{-i}=-2(\rho^i+m-i) \quad \textrm{on $V_{\rho-\mu_i}$}. \nonumber 
\end{equation}
So we have
$$
\widehat{C}_-(\phi\otimes \bar{\e}_k)=\widehat{C}_-(\sum_i p_{-i}(\bar{\e}_k)\phi) =\sum_i -2 w_{-i}p_{-i}(\bar{\e}_k)\phi.
$$
On the other hand, we have 
\begin{equation}
  \widehat{C}_-(\phi\otimes \bar{\e}_k) 
=2\sum_{s,t} \pi_{\rho}(e_{st})\phi \otimes \pi_{-\mu_m}(e_{t s})\bar{\e}_k 
=-2\sum_{l}p_{-i}(\bar{\e}_s)\pi_{\rho}(e_{sk})\phi
\nonumber
\end{equation}
for any $\phi$ in $V_{\rho}$. Then we conclude that 
\begin{equation}
w_{-i}p_{-i}(\bar{\e}_k)=\sum p_{-i}(\bar{\e}_l)\pi_{\rho}(e_{lk}) :V_{\rho}\to V_{\rho-\mu_i}. \nonumber
\end{equation}
Similarly, we consider the tensor product representation $(\pi_{\rho}\otimes \pi_{\mu_1},V_{\rho}\otimes \mathbb{C}^m)$ and the operator 
\begin{equation}
\widehat{C}_+:=\pi_{\rho}\otimes \pi_{\mu_1}(c_2)-\pi_{\rho}(c_2)\otimes \id -\id \otimes \pi_{\mu_1}(c_2). \label{eqn:3-14}
\end{equation}
In this case, we have
\begin{equation}
w_{+i}p_{+i}(\e_k)=-\sum p_{+i}(\e_l)\pi_{\rho}(e_{kl}) :V_{\rho}\to V_{\rho+\mu_i}. \nonumber
\end{equation}
%%%%%%%%%%%%%%%%%  lemma %%%%%%%%%%%%%%%%
\begin{lem}\label{lem:3-4}
Let $w_{+i}$ (resp. $w_{-i}$) be the $(1,0)$-conformal weight (resp. the $(0,1)$-conformal weight) associated to $\rho$. Then the Clifford homomorphisms $p_{\pm i}$ satisfy that
\begin{gather}
w_{+i}p_{+i}(\e_k)=-\sum p_{+i}(\e_l)\pi_{\rho}(e_{kl}) :V_{\rho}\to V_{\rho+\mu_i},\label{eqn:3-16}\\
w_{-i}p_{-i}(\bar{\e}_k)=\sum p_{-i}(\bar{\e}_l)\pi_{\rho}(e_{lk}) :V_{\rho}\to V_{\rho-\mu_i}. \label{eqn:3-17}
\end{gather}
\end{lem}
%%%%
From this lemma, we can relate the Clifford homomorphisms to the enveloping algebra. 
%%%%%%%%%%%%%%%%%%% theorem %%%%%%%%%%%%%%%%%%
\begin{thm}\label{thm:3-5}
Let $w_{\pm i}$ be the conformal weight, and $ \{ e_{kl}^q, \tilde{e}_{kl}^q \}_{k,l,q}$ be elements in $U(\gl)$ given in \eqref{eqn:2-1} and \eqref{eqn:2-9}. Then the Clifford homomorphisms satisfy that 
\begin{gather}
\sum_{1\le i\le m} w_{+i}^qp_{+i}(\e_k)^{\ast}p_{+i}(\e_l)=\pi_{\rho}(\tilde{e}_{kl}^q)  \in  \End (V_{\rho}), \label{eqn:3-18}\\
\sum_{1\le i\le m} w_{-i}^qp_{-i}(\bar{\e}_k)^{\ast}p_{-i}(\bar{\e}_l)=\pi_{\rho}(e_{kl}^q) \in  \End (V_{\rho}), \label{eqn:3-19}
\end{gather}
for $q$ in $\mathbb{Z}_{\ge 0}$ and $k,l=1,\cdots,m$. Furthermore, we take the trace of the above equations and have
\begin{gather}
\sum_{i,k}w_{+i}^qp_{+i}(\e_k)^{\ast}p_{+i}(\e_k)=\pi_{\rho}(\tilde{c}_q), \label{eqn:3-20} \\
\sum_{i,k} w_{-i}^qp_{-i}(\bar{\e}_k)^{\ast}p_{-i}(\bar{\e}_k)=\pi_{\rho}(c_q).\label{eqn:3-21}
\end{gather}
\end{thm}
%%%%%%%%% proof %%%%%%%%%%%%
\begin{proof}
We use the equation \eqref{eqn:3-16} $q$-times and have
\begin{equation}
w_{+i}^q p_{+i}(\e_l) =(-1)^q\sum_{i_1,\cdots,i_q} p_{+i}(\e_{i_1})\pi_{\rho}(e_{i_2i_1})\pi_{\rho}(e_{i_3i_2})\cdots \pi_{\rho}(e_{li_{q}}). \nonumber
\end{equation}
We multiply $p_{+i}(\e_k)^{\ast}$ on both sides and take the sum for $i$. Then we have from \eqref{eqn:3-10} 
\begin{equation}
\begin{split}
\sum_{i} w_{+i}^q p_{+i}(\e_k)^{\ast} p_{+i}(\e_l) &=(-1)^q\sum_{i_1,\cdots,i_{q-1}} \pi_{\rho}(e_{i_1k})\pi_{\rho}(e_{i_2i_1})\cdots \pi_{\rho}(e_{li_{q-1}}) \\ 
 &=\pi_{\rho}(\tilde{e}_{kl}^q).
\end{split}\nonumber
\end{equation}
\end{proof}
%%%%%%%%
This theorem means that $p_{+i}(\e_k)^{\ast}p_{+i}(\e_l)$ can be expressed as a linear combination of $\{\tilde{e}_{kl}^q\}_q$. We take a matrix expression of \eqref{eqn:3-18}, 
\begin{equation}
\begin{pmatrix}
w_{+1}^0  & w_{+2}^0 & \cdots & w_{+m}^0 \\
w_{+1}^1  & w_{+2}^1 & \cdots & w_{+m}^1 \\
 & \cdots \cdots &      &     \\
w_{+1}^{m-1} & w_{+2}^{m-1} &   \cdots & w_{+m}^{m-1}
\end{pmatrix}
\begin{pmatrix}
 p_{+1}(\e_k)^{\ast} p_{+1}(\e_l) \\
 p_{+2}(\e_k)^{\ast} p_{+2}(\e_l) \\
\cdots \\
p_{+m}(\e_k)^{\ast} p_{+m}(\e_l) 
\end{pmatrix} =
\begin{pmatrix}
\pi_{\rho}(\tilde{e}_{kl}^0) \\
\pi_{\rho}(\tilde{e}_{kl}^1) \\
\cdots \\
\pi_{\rho}(\tilde{e}_{kl}^{m-1})
\end{pmatrix}
 . \nonumber
\end{equation}
Because the $(1,0)$-conformal weights are different from each other, the above Vandermonde matrix has the inverse matrix whose $(i,j)$-th component is\begin{equation}
(-1)^{m-j}\frac{S_{m-j}(w_{+1},\cdots,\widehat{w_{+i}},\cdots,w_{+m})}{\prod_{j\neq i}(w_{+i}-w_{+j})}.\label{eqn:3-22}
\end{equation}
Here, $S_j(x_1,\cdots ,\hat{x}_i,\cdots,x_m)$ is the $j$-th fundamental symmetric polynomial of $\{ x_1,\cdots,x_m \}\setminus \{x_i\}$. 
%%%%%%%%%%%%  corollary %%%%%%%%%%%%%%%%
\begin{cor}\label{cor:3-6}
The Clifford homomorphisms are expressed as linear combinations of $\{e_{kl}^q\}_{q=0}^{m-1}$ and $\{\tilde{e}_{kl}^q\}_{q=0}^{m-1}$.
\begin{gather}
p_{+i}(\e_k)^{\ast} p_{+i}(\e_l)=\sum_{1\le j\le m}(-1)^{m-j}\frac{S_{m-j}(w_{+1},\cdots,\widehat{w_{+i}},\cdots,w_{+m})}{\prod_{j\neq i}(w_{+i}-w_{+j})}   \pi_{\rho}(\tilde{e}_{kl}^{j-1}), \label{eqn:3-23}\\
p_{-i}(\bar{\e}_k)^{\ast} p_{-i}(\bar{\e}_l)=\sum_{1\le j\le m}(-1)^{m-j}\frac{S_{m-j}(w_{-1},\cdots,\widehat{w_{-i}},\cdots,w_{-m})}{\prod_{j\neq i}(w_{-i}-w_{-j})}   \pi_{\rho}(e_{kl}^{j-1}). \label{eqn:3-24}
\end{gather}
We take the trace of these equations and have
\begin{gather}
\sum_k p_{+i}(\e_k)^{\ast} p_{+i}(\e_k)= \gamma_{+i}, \label{eqn:3-25}\\
\sum_k p_{-i}(\bar{\e}_k)^{\ast} p_{-i}(\bar{\e}_k)= \gamma_{-i}. \label{eqn:3-26}
\end{gather}
\end{cor}
%%%%%%
\begin{rem}
When $\rho\pm \mu_i$ is not dominant integral, we can show that $\gamma_{\pm i}$ is zero.
\end{rem}
From this corollary, our interest is to find out relations between $e_{kl}^q$ and $\tilde{e}_{kl}^q$, which are given in the next section.

Let us discuss relations between the orthogonal projection $\Pi_{\pm i}$ and the Clifford homomorphisms in the rest of this section. 
%%%%%%%%%%%%%%  proposition %%%%%%%%%%%%%%%%%%%%
\begin{prop}\label{prop:3-7}
We consider the orthogonal projection
$$
\Pi_{+i}:V_{\rho}\otimes \mathbb{C}^m\to V_{\rho+\mu_i}\subset V_{\rho}\otimes \mathbb{C}^m.
$$
This projection is given by
\begin{equation}
\Pi_{+i}(\phi\otimes u)=\sum_k p_{+i}(\e_k)^{\ast}p_{+i}(u)\phi \otimes \e_k
\in V_{\rho+\mu_i}\subset V_{\rho}\otimes \mathbb{C}^m \label{eqn:3-27}
\end{equation}
for $\phi$ in $V_{\rho}$ and $u$ in $\mathbb{C}^m$.
\end{prop}
%%%%%%%%%%% proof %%%%%%%%%%%%
\begin{proof}
There is an embedding $\Theta_{+i}$ from $V_{\rho+\mu_i}$ into $V_{\rho}\otimes \mathbb{C}^m$, 
$$
\Theta_{+i}:V_{\rho+\mu_i}\ni \phi\mapsto \Theta_{+i}(\phi)=\sum_k p_{+i}(\e_k)^{\ast}\phi\otimes \e_k \in V_{\rho}\otimes \mathbb{C}^m. 
$$
We know that $\Theta_{+i}$ does not depend the choice of unitary basis and is $U(m)$-equivariant. For $u=\sum u_l\e_l$ in $\mathbb{C}^m$, we have 
\begin{equation}
\begin{split}
\phi\otimes u &= \sum_l u_l \phi\otimes \e_l  \\
&=\sum_l u_l \sum_k \delta_{kl} \phi\otimes \e_k \\
 &=\sum_l u_l \sum_k \sum_{i} p_{+i}(\e_k)^{\ast}p_{+i}(\e_l)\phi\otimes \e_k \\ &=\sum_i \sum_k p_{+i}(\e_k)^{\ast}p_{+i}(u) \phi\otimes \e_k=\sum_{i} \Theta_{+i}(p_{+i}(u)\phi).
\end{split}\nonumber
\end{equation}
Then we have proved the proposition. 
\end{proof}
%%%%%%%%%%% lemma %%%%%%%%%%%%%%%%%%%
\begin{lem}\label{lem:3-8}
Suppose that $\rho+\mu_i$ is dominant integral, and we have
\begin{equation}
\sum_k p_{+i}(\e_k)p_{+i}(\e_k)^{\ast}=\id :V_{\rho+\mu_i}\to V_{\rho+\mu_i}. \label{eqn:3-28}
\end{equation}
\end{lem}
%%%%%%%%% proof %%%%%%%%%%%%%%
\begin{proof}
We can show that $\sum_k p_{+i}(\e_k)p_{+i}(\e_k)^{\ast}$ is a $U(m)$-equivariant endomorphism of $V_{\rho+\mu_i}$, and hence, is constant. We calculate the constant. 
\begin{equation}
\begin{split}
 &(p_{+i}(u)\psi,p_{+i}(v)\phi)  \\
=& \sum_{k,l} (p_{+i}(\e_k)^{\ast}p_{+i}(u)\psi\otimes \e_k,  p_{+i}(\e_l)^{\ast}p_{+i}(v)\psi\otimes \e_l) \\
=& \sum_{k,l} (p_{+i}(\e_l)p_{+i}(\e_k)^{\ast}p_{+i}(u)\psi, p_{+i}(v)\phi)\delta_{kl} \\
=&( \sum_{k} p_{+i}(\e_k)p_{+i}(\e_k)^{\ast}p_{+i}(u)\psi, p_{+i}(v)\phi)
  \\ 
=&\sum_{k} p_{+i}(\e_k)p_{+i}(\e_k)^{\ast} (p_{+i}(u)\psi, p_{+i}(v)\phi)
\end{split} \nonumber
\end{equation}
for any $\phi, \psi$ in $V_{\rho+\mu_i}$ and $u,v$ in $\mathbb{C}^m$. So we have $\sum_{k} p_{+i}(\e_k)p_{+i}(\e_k)^{\ast}=1$. 
\end{proof}
%%%%%%%

When $\rho+\mu_i$ is dominant integral, namely $V_{\rho+\mu_i}\neq \{0\}$, we have the Clifford homomorphism $p^{\rho}_{\rho+\mu_i}(u)^{\ast}=p_{+i}(u)^{\ast}:V_{\rho+\mu_i}\to V_{\rho}$. But, there is another Clifford homomorphism $p^{\rho+\mu_i}_{\rho}(\bar{u}):V_{\rho+\mu_i}\to V_{\rho}$. Then we show that 
%%%%%%%%%%%%%%%%% lemma %%%%%%%%%%%%
\begin{prop}\label{prop:3-9}
The Clifford homomorphism $p_{\rho+\mu_i}^{\rho}(u)^{\ast}:V_{\rho+\mu_i}\to V_{\rho}$ and $p^{\rho+\mu_i}_{\rho}(\bar{u}):V_{\rho+\mu_i}\to V_{\rho}$ are related as follows:
\begin{equation}
p^{\rho+\mu_i}_{\rho}(\bar{u})=\frac{1}{\sqrt{\gamma_{+i}}}p_{\rho+\mu_i}^{\rho}(u)^{\ast}:V_{\rho+\mu_i}\to V_{\rho}. \label{eqn:3-29}
\end{equation}
Similarly, we have 
\begin{equation}
p^{\rho-\mu_i}_{\rho}(u)=\frac{1}{\sqrt{\gamma_{-i}}}p_{\rho-\mu_i}^{\rho}(\bar{u})^{\ast}:V_{\rho-\mu_i}\to V_{\rho}. \label{eqn:3-30}
\end{equation}
Here, there is ambiguity of complex number with norm $1$ on each equation.
\end{prop}
%%%%%%%%%%%%%%%
\begin{proof}
We have an embedding from $V_{\rho+\mu_i}$ to $V_{\rho}\otimes \mathbb{C}^m$, 
\begin{equation}
\Theta_{+i}':V_{\rho+\mu_i}\ni\phi\mapsto \Theta'_{+i}(\phi)=\sum p^{\rho+\mu_i}_{\rho}(\bar{\e}_k)\phi\otimes \e_k\in V_{\rho}\otimes \mathbb{C}^m. \nonumber
\end{equation}
It is from Schur's lemma that there is a complex number $a$ such that $\Theta_{+i}'=a\Theta_{+i}$, that is,
$$
\sum p^{\rho+\mu_i}_{\rho}(\bar{\e}_k)\phi\otimes \e_k=a\sum p_{\rho+\mu_i}^{\rho}(\e_k)^{\ast}\phi\otimes \e_k
$$
for any $\phi$ in $V_{\rho+\mu_i}$. Then we have $p^{\rho+\mu_i}_{\rho}(\bar{u})=ap_{\rho+\mu_i}^{\rho}(u)^{\ast}$. Furthermore, we know that 
\begin{equation}
\gamma_{+i}=\sum_k p_{\rho+\mu_i}^{\rho}(\e_k)^{\ast}p_{\rho+\mu_i}^{\rho}(\e_k)= |a|^{-2} \sum p^{\rho+\mu_i}_{\rho}(\bar{\e}_k)p^{\rho+\mu_i}_{\rho}(\bar{\e}_k)^{\ast}=|a|^{-2}. \nonumber
\end{equation}
Thus, it holds that $p^{\rho+\mu_i}_{\rho}(\bar{u})=\frac{ e^{\sqrt{-1}\theta}}{\sqrt{\gamma_{+i}}}p_{\rho+\mu_i}^{\rho}(u)^{\ast}$. The ambiguity $e^{\sqrt{-1}\theta}$is due to the discussion in Remark \ref{rem:3-2}.
\end{proof}

%%%%%%%%%%%%%%%%%%%%%%%%%%%%%%%
%%            4               %%
%%%%%%%%%%%%%%%%%%%%%%%%%%%%%%%
\section{Algebraic relations for Clifford homomorphisms}\label{sec:4}
In this section, we  have further relations for the Clifford homomorphisms. From theorem \ref{thm:3-5} or corollary \ref{cor:3-6}, our problem reduces to construct relations between $\{e_{kl}^q\}_q$ and $\{\tilde{e}_{kl}^q\}_q$.

It follows from lemma \ref{lem:2-1} and \ref{lem:2-3} that 
\begin{equation}
\begin{split}
\sum e_{il}^qe_{ki}=e_{kl}^{q+1}-me_{kl}^q+\delta_{kl}c_q, \\
\sum \tilde{e}_{il}^q\tilde{e}_{ki}=\tilde{e}_{kl}^{q+1}-m\tilde{e}_{kl}^q+\delta_{kl}\tilde{c}_q.
\end{split}  \label{eqn:4-1}
\end{equation}
The next lemma allows us to represent $\tilde{e}_{kl}^q$ as a linear combination of $\{e_{lk}^p\}_p$ whose coefficients are polynomials of the Casimir elements. %%%%%%%%%%%%%%%%%%  lem %%%%%%%%%%%%%%%%
\begin{lem}\label{lem:4-1}
Let $e_{kl}^q$ and $\tilde{e}_{lk}^q$ be the elements in $U(\gl)$ defined in \eqref{eqn:2-1} and \eqref{eqn:2-9}. There exists a set of polynomials $\{b_{q,p}(c)| 0\le p\le q\}$ of Casimir elements $c=(c_0,c_1,\cdots )$ such that
\begin{equation}
\tilde{e}_{kl}^q=\sum_{p=0}^q b_{q,p}(c) e_{lk}^p \label{eqn:4-2}
\end{equation}
 for any $k,l$. 
\end{lem}
%%%%%%
\begin{proof}
We prove the lemma by an induction for $q$. When $q$ is $0$ and $1$, the lemma is trivial. We assume that the lemma holds for non-negative integers $\le q$. Then, for any $k,l$, 
\begin{equation}
\begin{split}
 \tilde{e}_{kl}^{q+1}=&\sum_i \tilde{e}_{ki}^q\tilde{e}_{il}=-\sum_i\tilde{e}_{ki}^qe_{li}=-\sum_i \sum_{p=0}^qb_{q,p}e_{ik}^p e_{li} \\    
 =&- \sum_{p=0}^q b_{q,p}( e_{lk}^{p+1}-me_{lk}^p+\delta_{lk} c_p) \quad \textrm{(by \eqref{eqn:4-1})}\\
    =&-b_{q,q}e_{lk}^{q+1}+\sum_{p=1}^q(m b_{q,p}-b_{q,p-1})e_{lk}^p+(m b_{q,0}-\sum_{p=0}^q b_{q,p}c_p)\delta_{lk} \\
    =&\sum_{p=0}^{q+1} b_{q+1,p}e_{lk}^p.
\end{split}\nonumber
\end{equation}
Thus there exists a set of polynomials $\{b_{q,p}\}_{q,p}$ of the Casimir elements, which are defined inductively. 
\end{proof}
In the above proof, we have a recursion formula for $\{b_{q,p}\}_{q,p}$. Since it is a little difficult to solve the recursion formula, we consider the equation 
\begin{equation}
\sum_{p=0}^q 
 \binom{q}{p} (-m)^{q-p}
\tilde{e}_{kl}^p=\sum_{p=0}^q a_{q,p}(c) e_{lk}^p, \label{eqn:4-3}
\end{equation}
and a recursion formula for $\{a_{q,p}\}_{q,p}$. 
%%%%%%%%%%%%%%%%% lemma %%%%%%%%%%%%%%%%
\begin{lem}\label{lem:4-2}
Let $a_{q,p}(c)$ be a polynomial of the Casimir elements $c=(c_0,c_1,\cdots)$ satisfying \eqref{eqn:4-3}. Then $\{a_{q,p}\}_{0\le p\le q}$ are given inductively by the recursion formula 
\begin{equation}
a_{q,p}=(-1)^pa_{q-p,0}, \quad a_{q,q}=(-1)^q, \quad a_{q,0}=-\sum_{p=0}^{q-1}a_{p,0}x_{q-p}, \label{eqn:4-4}
\end{equation}
where we set $x_p:=(-1)^{p-1}c_{p-1}$  for $p=1,2,\cdots$.
\end{lem}
%%%%%%%%%%%%%%%%
\begin{proof}
\begin{equation}
\begin{split}
&\sum_{p=0}^{q+1} 
\binom{q+1}{p}
 (-m)^{q+1-p}
\tilde{e}_{kl}^p 
=\sum_{p=0}^{q+1} 
\left( \binom{q}{p}+ \binom{q}{p-1}\right)
\tilde{e}_{kl}^p(-m)^{q+1-p} \\
=&-m\sum_{p=0}^q a_{q,p}e_{lk}^p+\sum_i \sum_{p=1}^{q+1} \binom{q}{p-1}  \tilde{e}_{ki}^{p-1}\tilde{e}_{il}^1 (-m)^{q+1-p} \\
=&-m\sum_{p=0}^q a_{q,p}e_{lk}^p+\sum_i \sum_{p=0}^{q} a_{q,p}e_{ik}^p(-e_{li})\\
=&-m\sum_{p=0}^q a_{q,p}e_{lk}^p+\sum_{p=0}^qa_{q,p}(-e_{lk}^{p+1}+me_{lk}^p-\delta_{lk}c_p) \\
=&-\sum_{p=1}^{q+1}a_{q,p-1}e_{lk}^p-(\sum_{p=0}^q a_{q,p}c_p)\delta_{lk}
=\sum_{p=0}^{q+1}a_{q+1,p}e_{lk}^p.
\end{split}\nonumber
\end{equation}
Then we have a recursion formula of $\{a_{q,p}\}_{0\le p\le q}$,
\begin{equation}
a_{q,p}=-a_{q-1,p-1}, \quad a_{0,0}=1, \quad a_{q,0}=-\sum_{p=0}^{q-1}a_{q-1,p}c_p.\nonumber
\end{equation}
This formula is equivalent to \eqref{eqn:4-4}.
\end{proof}
%%%%%%%%%
To give the solution of \eqref{eqn:4-4}, we prepare a generating function of valuable $z$ whose coefficients are polynomials of $x=(x_1,x_2,\cdots)$. 
We define a generating function $K(z)$ by
\begin{equation}
K(z):=\sum_{n=0}^{\infty} K_{n}(x)z^n=\frac{1}{1+x_1z^1+x_2z^2+\cdots}.\label{eqn:4-5}
\end{equation}
Here, each coefficient of $z$ is 
\begin{equation}
K_n(x)=\sum_{i_1+2i_2+\cdots +ni_n=n}\frac{(i_1+i_2+\cdots+i_n)!}{i_1!i_2!\cdots i_n!}(-x_1)^{i_1}\cdots (-x_n)^{i_n} \label{eqn:4-6}
\end{equation}
(see \cite{MS}). We know that the first few terms are  
\begin{equation}
K_0(x) =1, \quad K_1(x) =-x_1, \quad K_2(x)=x_1^2-x_2, \quad 
K_3(x) =-x_1^3+2x_1x_2-x_3,
\nonumber
\end{equation}
and $K_n(x)$ satisfies 
\begin{equation}
K_n(x_1,-x_2,x_3,-x_4\cdots)=(-1)^n K_n(-x). \label{eqn:4-7}
\end{equation}
Let us give a recursion formula of $\{K_{n}(x)\}_n$. Because of the definition of $K(z)$, we have 
\begin{equation}
1 =(\sum_{i=0}^{\infty} x_iz^i)(\sum_{j=0}^{\infty}K_j(x)z^j)
  =\sum_{q=0}^{\infty} (\sum_{p=0}^{q} K_{p}(x)x_{q-p}) z^q.
\nonumber
\end{equation}
Compare the coefficients of $z^n$, and we have
\begin{equation}
K_0(x) =1, \quad K_{q}(x)=-\sum_{p=0}^{q-1} K_{p}(x)x_{q-p}. \label{eqn:4-8}
\end{equation}
This is just the recursion formula \eqref{eqn:4-4} of $\{a_{p,0}\}_{p\ge 0}$. Then we conclude that
\begin{equation}
a_{q,p}=(-1)^pa_{q-p,0}=(-1)^pK_{q-p}(c_0,-c_1,c_2,\cdots)=(-1)^qK_{q-p}(-c). 
\nonumber %\label{eqn:4-9}
\end{equation}
Thus, we have the following theorem.
%%%%%%%%%%%%%%%%%%%%%%%%%% thm %%%%%%%%%%%%%%%%%%
\begin{thm}\label{thm:4-3}
Let $\{e_{kl}^p\}_p$ and $\{\tilde{e}_{kl}^p\}_p$ be the elements of $U(\gl)$ given in \eqref{eqn:2-1} and \eqref{eqn:2-9}. Then it holds that
\begin{equation}
\sum_{p=0}^q 
\binom{q}{p} (-m)^{q-p}
\tilde{e}_{kl}^p=(-1)^q\sum_{p=0}^q K_{q-p}(-c) e_{lk}^p, \label{eqn:4-10}
\end{equation}
where 
\begin{equation}
K_n(-c)=\sum_{i_1+2i_2+\cdots +ni_n=n}\frac{(i_1+i_2+\cdots+i_n)!}{i_1!i_2!\cdots i_n!}c_0^{i_1}\cdots c_{n-1}^{i_n}. \label{eqn:4-11}
\end{equation}
We take the involution of the above relation and have
\begin{equation}
\sum_{p=0}^q 
\binom{q}{p} (-m)^{q-p}
e_{kl}^p=(-1)^q\sum_{p=0}^q K_{q-p}(-\tilde{c}) \tilde{e}_{lk}^p. \label{eqn:4-12}
\end{equation}
\end{thm}
%%%%%%%%%
Taking the trace of \eqref{eqn:4-10}, we have a relation between the Casimir elements $\{c_p\}_p$ and $\{\tilde{c}_p\}_p$. 
%%%%%%%%%%%%%%%%%%% cor %%%%%%%%%%%%%%%%%%%
\begin{cor}\label{cor:4-4}
The Casimir elements $\{c_p \}_p$ and $\{ \tilde{c}_p\}_p$ are related to each other as follows:
\begin{equation}
\begin{split}
\sum_{p=0}^q 
\binom{q}{p} (-m)^{q-p}
\tilde{c}_p=(-1)^q K_{q+1}(-c), \\
 \sum_{p=0}^q \binom{q}{p} (-m)^{q-p}
c_p=(-1)^{q}K_{q+1}(-\tilde{c}).
\end{split} \label{eqn:4-13}
\end{equation}
\end{cor}
%%%%%%%%%
The above theorem induces relations for the Clifford homomorphisms.
%%%%%%%%%%%  cor %%%%%%%%%%%%%%%%%%%%%
\begin{cor}\label{cor:4-5}
The Clifford homomorphisms $p_{\pm i}$ and $p_{\pm i}^{\ast}$ satisfy that
\begin{equation}
\begin{split}
 &\sum_{1\le i\le m}(w_{+i}-m)^q p_{+i}(\e_k)^{\ast}p_{+i}(\e_l) 
  \\=&(-1)^q \sum_{1\le i \le m} (\sum_{p=0}^q \pi_{\rho}(K_{q-p}(-c))w_{-i}^p) p_{-i}(\bar{\e}_l)^{\ast}p_{-i}(\bar{\e}_k), \\
 &\sum_{1\le i\le m}(w_{-i}-m)^q p_{-i}(\bar{\e}_k)^{\ast}p_{-i}(\bar{\e}_l) 
  \\=&(-1)^q \sum_{1\le i \le m} (\sum_{p=0}^q \pi_{\rho}(K_{q-p}(-\tilde{c}))w_{+i}^p) p_{+i}(\e_l)^{\ast}p_{+i}(\e_k)  
\end{split}\label{eqn:4-14}
\end{equation}
on irreducible $U(m)$-module $V_{\rho}$ with highest weight $\rho$.
\end{cor}
%%%%%%%%%
%%%%%%%
\begin{rem}\label{rem:4-1}
This corollary says that $p_{+i}(\e_k)^{\ast}p_{+i}(\e_l)$ is a linear combination of $\{p_{-i}(\bar{\e}_l)^{\ast}p_{-i}(\bar{\e}_k)\}_{i=1}^m$. Moreover, we remark that, if there are $N$ irreducible components in $V_{\rho}\otimes \mathbb{C}^m$, then we have $N$ components in $V_{\rho}\otimes \overline{\mathbb{C}^m}$, and vice versa. So, if there are $2N$ Clifford homomorphisms on $V_{\rho}$, then at most $N$ relations make sense in \eqref{eqn:4-14}. 
\end{rem}
%%%%%%%%%%%%%%%
It follows from \eqref{eqn:4-13} that $\tilde{c}_q$ is a polynomial of $\{c_p\}_p$. Similarly, we can express $\tilde{e}_{kl}^q$ as a linear combination of $\{e_{lk}^p\}_p$ whose coefficients are polynomials of $\{c_p\}_p$. The following corollary is the answer to the recursion formula for $\{b_{q,p}\}_{q,p}$ in lemma \ref{lem:4-1}. 
%%%%%%%%%%%%%%  cor %%%%%%%%%%%%%%%%%
\begin{cor}\label{cor:4-6}
We can express $\tilde{e}_{kl}^q$ and $\tilde{c}_q$ with $\{e_{lk}^p\}_p$ and $\{c_p\}_p$ as follows:
\begin{gather}
\tilde{e}_{kl}^q=(-1)^q \sum_{p=0}^q \{ \sum_{s=p}^q \binom{q}{s}  (-m)^{q-s} K_{s-p}(-c) \} e_{lk}^p, \label{eqn:4-15}\\
\tilde{c}_q=(-1)^q\sum_{p=0}^{q} \binom{q}{p} (-m)^{q-p} K_{p+1}(-c). \label{eqn:4-16}
\end{gather}
\end{cor}
%%%%%%%
%%%%%%%%%%%%%%%%%%%%%%%%%%%%%%%%
%%            5               %%
%%%%%%%%%%%%%%%%%%%%%%%%%%%%%%%%
\section{Example: spinors and Clifford multiplications} \label{sec:5}
We shall discuss spinors as an example and show that the Clifford homomorphisms on spinor space are the usual Clifford multiplications. Let $\mathbb{R}^{2m}\simeq \mathbb{C}^m$ be the Euclidean space with standard complex structure $J$ and $\Lambda^{0,p}$ be the space of $(0,p)$-forms. This space $\Lambda^{0,p}=\Lambda^{p}\overline{(\mathbb{C}^m)^{\ast}}$ is irreducible $U(m)$-module with highest weight $\rho=(1_p,0_{m-p})$, and is isomorphic to $\Lambda^p(\mathbb{C}^m)$. The direct sum $\sum \Lambda^{0,p}$ is known as the spinor space associated to $\mathbb{C}^{m}$. 

We have four Clifford homomorphisms from $V_{\rho}=\Lambda^{0,p}$ to other $U(m)$-modules,
\begin{equation}
\begin{split}
p_{+1}(\e_k) &:\Lambda^{0,p}\to V_{\rho+\mu_1}, \\
p_{+(p+1)}(\e_k)&:\Lambda^{0,p}\to V_{\rho+\mu_{p+1}}=\Lambda^{0,p+1}, \\
p_{-m}(\bar{\e}_k)&:\Lambda^{0,p}\to V_{\rho-\mu_m}, \\
p_{-p}(\bar{\e}_k)&:\Lambda^{0,p}\to V_{\rho-\mu_p}=\Lambda^{0,p-1}.
\end{split}\nonumber % \label{eqn:5-1}
\end{equation}
The conformal weights and $\gamma_{\pm i}$ associated to each Clifford homomorphism are 
\begin{equation}
\begin{split}
w_{+1}=-1,    &\quad   \gamma_{+1}=\frac{p(m+1)}{p+1}, \\
w_{+(p+1)}=p,  &\quad  \gamma_{+(p+1)}=\frac{m-p}{p+1}, \\
w_{-m}=0 ,     &\quad  \gamma_{-m}=\frac{(m+1)(m-p)}{m-p+1},\\
w_{-p}=m-p+1,   &\quad \gamma_{-p}=\frac{p}{m-p+1}.
\end{split}\nonumber % \label{eqn:5-2}
\end{equation}
From theorem \ref{thm:3-5}, it holds that
\begin{gather}
p_{+1}(\e_k)^{\ast}p_{+1}(\e_l)+p_{p+1}(\e_k)^{\ast}p_{p+1}(\e_l)
=\delta_{kl},\nonumber \\ 
p_{-m}(\bar{\e}_l)^{\ast}p_{-m}(\bar{\e}_k)+p_{-p}(\bar{\e}_l)^{\ast}p_{-p}(\bar{\e}_k)=\delta_{kl}, \nonumber\\
-p_{+1}(\e_k)^{\ast}p_{+1}(\e_l)+pp_{+(p+1)}(\e_k)^{\ast}p_{+(p+1)}(\e_l)=-\pi_{\rho}(e_{lk}), \nonumber\\
(m-p+1)p_{-p}(\bar{\e}_k)^{\ast}p_{-p}(\bar{\e}_l)=\pi_{\rho}(e_{kl}).\nonumber
\end{gather}
Here, we remark that other relations reduce to these equations. From the above 
equations, we have
\begin{equation}
(p+1)p_{+(p+1)}(\e_k)^{\ast}p_{+(p+1)}(\e_l)+(m-p+1)p_{-p}(\bar{\e}_l)^{\ast}p_{-p}(\bar{\e}_k)=\delta_{kl}. \label{eqn:5-3}
\end{equation}
We shall prove that this relation coincides with usual Clifford relation.
We set, for any $p$,
\begin{equation}
\begin{split}
-\bar{\e}_k\cdot:=i(\bar{\e}_k)=\sqrt{m-p+1}\: p_{-p}(\bar{\e}_k) :\Lambda^{0,p}\to \Lambda^{0,p-1}, \\
\e_k\cdot:=\e_k{}_{\wedge}=\sqrt{p+1}\:p_{+(p+1)}(\e_k):\Lambda^{0,p}\to \Lambda^{0,p+1}.
\end{split}\label{eqn:5-4}
\end{equation}
It follows from \eqref{eqn:3-29} and \eqref{eqn:3-30} that we have
$$
(\bar{\e}_k)^{\ast}=-\e_k, \quad (\e_k)^{\ast}=-\bar{\e}_k.
$$
Thus, the equation \eqref{eqn:5-3} means usual Clifford relation $\e_k\bar{\e}_l+\bar{\e}_l\e_k=-\delta_{kl}$. Furthermore, we have 
\begin{equation}
-\e_k\bar{\e}_l=\pi_{\rho}(e_{kl}) \label{eqn:5-5}
\end{equation}
on each $U(m)$-module $\Lambda^{0,p}$. 
%%%%%%%%%%%%%%%%% rem %%%%%%%%%%%%%%
\begin{rem}
How can we get other Clifford relations $\e_k\e_l+\e_l\e_k=0$ and $\bar{\e}_k\bar{\e}_l+\bar{\e}_l\bar{\e}_k=0$ ? We consider the orthogonal projection from $\Lambda^{0,p}\otimes S^2(\mathbb{C}^m)$ onto $\Lambda^{0,p+2}$. By using the tensor product decomposition rule, we have 
$$
\Lambda^{0,p}\otimes S^2(\mathbb{C}^m)=V_{\rho+2\mu_1}\oplus V_{\rho+\mu_1+\mu_{p+1}}.
$$
So $\Lambda^{0,p+2}$ does not appear as an irreducible component of $\Lambda^{0,p}\otimes S^2(\mathbb{C}^m)$. Thus, we have the relation $\e_k\e_l+\e_l\e_k=0$ for $k,l=1,\cdots,m$. 
\end{rem}
%%%%%%%%%%%%%%%

Let us show an explicit formula of the orthogonal projection $\Pi_{\pm i}$, which is useful to construct (K\"ahlerian) twistor spinors. Considering the projection formula \eqref{eqn:3-27}, we have 
\begin{equation}
\begin{split}
 &\phi\otimes \e_k \\
=&\sum_l p_{+(p+1)}(\e_l)^{\ast}p_{+(p+1)}(\e_k)\phi\otimes \e_l+ p_{+1}(\e_l)^{\ast}p_{+1}(\e_k)\phi \otimes \e_l\\
=&-\frac{1}{p+1}\sum_{l} \bar{\e}_l\cdot\e_k\cdot\phi\otimes \e_l +(\phi\otimes \e_k+\frac{1}{p+1}\sum_l \bar{\e}_l\cdot\e_k\cdot\psi\otimes \e_l)
\end{split}\label{eqn:5-6}
\end{equation}
for $\phi$ in $\Lambda^{0,p}$. Similarly, we have
\begin{equation}
\begin{split}
\phi\otimes \bar{\e}_k &=
 -\frac{1}{m-p+1}\sum_l\e_l\cdot\bar{\e}_k\cdot\phi \otimes \bar{\e}_l+(\phi\otimes \bar{\e}_k+\frac{1}{m-p+1}\sum_l \e_l\cdot\bar{\e}_k\cdot\phi \otimes \bar{\e}_l).
\end{split}\label{eqn:5-7}
\end{equation}
%%%%%%%%%%%%%%%%%%%%%%%%%%%%%%%%
%%            6              %%
%%%%%%%%%%%%%%%%%%%%%%%%%%%%%%%%
\section{The K\"ahlerian gradients and their conformal covariance} \label{sec:6} In this section, we shall define the K\"ahlerian gradients on almost Hermitian manifolds, and show their conformal covariance. 

Let $M$ be a real $2m$-dimensional almost Hermitian manifold with almost complex structure $J$ and Hermitian metric $g$. Here, a Riemannian metric $g$ is called Hermitian metric if $g$ is compatible with $J$. The complexification of the tangent bundle splits into the direct sum of the $(1,0)$-tangent bundle $T^{1,0}(M)$ and the $(0,1)$-tangent bundle $T^{0,1}(M)$ with respect to $J$, where each bundle is equipped with Hermitian metric. Since the almost Hermitian structure gives a reduction to $U(m)$ of the structure group of frame bundle on $M$, we have a principal bundle $\mathbf{U}(M)$ with structure group $U(m)$ which is a bundle of the unitary frames on $M$. 

Let $(\pi_{\rho},V_{\rho})$ be an irreducible unitary representation of $U(m)$ with highest weight $\rho$ and $\mathbf{S}_{\rho}$ be the associated bundle $\mathbf{U}(M)\times_{\rho}V_{\rho}$ whose fiber metric is induced from the inner product on $V_{\rho}$. We fix a connection $\omega$ on $\mathbf{U}(M)$, namely, a linear connection compatible with $g$ and $J$. This connection gives a covariant derivative $\nabla$ preserving the metric on $\mathbf{S}_{\rho}$,
\begin{equation}
\nabla:\Gamma (\mathbf{S}_{\rho})\to \Gamma (\mathbf{S}_{\rho}\otimes T^{\ast}(M))=\Gamma (\mathbf{S}_{\rho}\otimes (T^{\ast}(M)\otimes \mathbb{C})). \label{eqn:6-1}
\end{equation}
The complexified cotangent bundle $T^{\ast}(M)\otimes \mathbb{C}$ splits as $\Lambda^{1,0}(M)\oplus \Lambda^{0,1}(M)$ and hence the connection $\nabla$ also splits as $\nabla^{1,0}+\nabla^{0,1}$. Regard $\Lambda^{1,0}(M)$ and $\Lambda^{0,1}(M)$ as $T^{0,1}(M)$ and $T^{1,0}(M)$ by $g$, and we have descriptions of $\nabla^{1,0}$ and $\nabla^{0,1}$ as follows: 
\begin{equation}
\begin{split}
\nabla^{1,0}\phi=\sum_{1\le k\le m} \nabla_{\e_k}\phi \otimes \bar{\e}_k \in \Gamma(\mathbf{S}_{\rho}\otimes T^{0,1}(M)), \\
\nabla^{0,1}\phi=\sum_{1\le k\le m} \nabla_{\bar{\e}_k}\phi \otimes \e_k \in \Gamma(\mathbf{S}_{\rho}\otimes T^{1,0}(M))
\end{split}\label{eqn:6-2}
\end{equation}
for $\phi$ in $\Gamma (\mathbf{S}_{\rho})$. Here, $\{e_k,Je_k\}_{k=1}^m$ is a local orthonormal frame of $T(M)$, and $\{\e_k, \bar{\e}_k\}_{k=1}^m$ are local unitary frames of $T^{1,0}(M)$ and $T^{0,1}(M)$ given as
\begin{equation}
\e_k=\frac{1}{\sqrt{2}}(e_k-\sqrt{-1}Je_k), \quad \bar{\e}_k=\frac{1}{\sqrt{2}}(e_k+\sqrt{-1}Je_k). \label{eqn:6-3}
\end{equation}

We have already seen the decomposition $V_{\rho}\otimes \mathbb{C}^m=\sum V_{\rho+\mu_i}$ and $V_{\rho}\otimes \overline{\mathbb{C}^m}=\sum V_{\rho-\mu_i}$ in section \ref{sec:3}. Hence, the tensor product bundle $\mathbf{S}_{\rho}\otimes T^{1,0}(M)$ and $\mathbf{S}_{\rho}\otimes T^{0,1}(M)$ are decomposed as 
\begin{equation}
\begin{split}
\mathbf{S}_{\rho}\otimes \Lambda^{0,1}(M)=\mathbf{S}_{\rho}\otimes T^{1,0}(M)=\sum_i \mathbf{S}_{\rho+\mu_i},\\
\mathbf{S}_{\rho}\otimes \Lambda^{1,0}(M)=\mathbf{S}_{\rho}\otimes T^{0,1}(M)=\sum_i \mathbf{S}_{\rho-\mu_i}.
\end{split}\label{eqn:6-4}
\end{equation}
Here, when $\rho\pm \mu_i$ is not dominant integral, we put $\mathbf{S}_{\rho \pm \mu_i}:=M\times \{0\}$. 
%%%%%%%%%%%%%%%%%%%  definition %%%%%%%%%%%%%%%%%%%%%%%
\begin{defini}\label{def:6-1}
Let $\mathbf{S}_{\rho}$, $\mathbf{S}_{\rho \pm \mu_i}$, $\nabla^{1,0}$, and $\nabla^{0,1}$ be as above. We define a set of first order differential operators $\{D_{\pm i}\}_{i=1}^m$ to be 
\begin{equation}
\begin{split}
D_{+i}:\Gamma (\mathbf{S}_{\rho})\xrightarrow{\nabla^{0,1}} \Gamma (\mathbf{S}_{\rho}\otimes \Lambda^{0,1}(M))\xrightarrow{\simeq}\Gamma (\mathbf{S}_{\rho}\otimes T^{1,0}(M))  \xrightarrow{\Pi_{+i}} \Gamma (\mathbf{S}_{\rho+\mu_i}), \\
D_{-i}:\Gamma (\mathbf{S}_{\rho})\xrightarrow{\nabla^{1,0}} \Gamma (\mathbf{S}_{\rho}\otimes \Lambda^{1,0}(M))\xrightarrow{\simeq}\Gamma (\mathbf{S}_{\rho}\otimes T^{0,1}(M)) \xrightarrow{\Pi_{-i}} \Gamma (\mathbf{S}_{\rho-\mu_i}), 
\end{split}\label{eqn:6-5}
\end{equation}
where $\Pi_{\pm i}$ is the orthogonal projection to $\mathbf{S}_{\rho\pm \mu_i}$ defined fiberwise. We call these operators \textit{the K\"ahlerian gradients} associated to $\rho$.
\end{defini}
%%%%%%%%%%
%%%%%%%%%%%  rem %%%%%%%%%%%%%%%
\begin{rem}
If we consider a diffeomorphism $\phi:M \mapsto  M'$ preserving almost Hermitian structure and connection, then $\phi(D_{\pm i})=D_{\pm i}'$. Thus, the K\"ahlerian gradients are the $U(m)$-invariant first order differential operators.
\end{rem}
%%%%%%%%%%%%

To give a formula of $D_{\pm i}$, we lift the Clifford homomorphisms to bundle homomorphisms. Since $T^{1,0}(M)$ and $\mathbf{S}_{\rho}$ are vector bundle associated to $\mathbf{U}(M)$, we define the Clifford (bundle) homomorphisms from $\mathbf{S}_{\rho}$ to $\mathbf{S}_{\rho+\mu_i}$ by 
$$
(T^{1,0}(M))_x\times (\mathbf{S}_{\rho})_x \ni [p, u]\times [p,\phi]\mapsto [p, p_{+i}(u)\phi] \in (\mathbf{S}_{\rho+\mu_i})_x,
$$
where $x$ is in $M$ and $p$ is in $\mathbf{U}(M)$. We know from proposition \ref{prop:3-2} that this map is well-defined. Thus, we have bundle homomorphisms $p_{+i}(X)\in \Gamma (\Hom (\mathbf{S}_{\rho},\mathbf{S}_{\rho+\mu_i}))$ and $p_{-i}(\overline{X})\in \Gamma (\Hom (\mathbf{S}_{\rho},\mathbf{S}_{\rho-\mu_i}))$ for each $X$ in $\Gamma (T^{1,0}(M))$. 
%%%%%%%%%%%%  proposition %%%%%%%%%%%
\begin{prop}\label{prop:6-2}
Let $D_{\pm i}$ be the K\"ahlerian gradient and $p_{\pm i}$ be the Clifford homomorphism. Then 
\begin{equation}
\begin{split}
D_{+i}=\sum p_{+i}(\e_k)\nabla_{\bar{\e}_k} :\Gamma (\mathbf{S}_{\rho})\to \Gamma (\mathbf{S}_{\rho+\mu_i}), \\
D_{-i}=\sum p_{-i}(\bar{\e}_k)\nabla_{\e_k} :\Gamma (\mathbf{S}_{\rho})\to \Gamma (\mathbf{S}_{\rho-\mu_i}).
\end{split}\label{eqn:6-6}
\end{equation}
\end{prop}
%%%%%%%%%

The following results in the theory of the Dirac operator can be generalized to the ones for $D_{\pm i}$:
\begin{enumerate}
	\item The Dirac operator $D$ satisfies $[D,f]=D\circ f-fD=\mathrm{grad}f\cdot$ for $f$ in $C^{\infty}(M)$. Here, $\mathrm{grad}f\cdot$ is the Clifford multiplication of gradient vector field of $f$.
 \item The connection on spinor bundle $\mathbf{S}(M)$ is compatible with the Levi-Civita connection, that is, $\nabla_V (W\cdot \phi)=W\cdot (\nabla_V \phi)+(\nabla_VW)\cdot \phi$ for $\phi$ in $\Gamma (\mathbf{S}(M))$ and $V,W$ in $\Gamma (T(M))$.
 \item The Dirac operator has a conformal covariance. 
\end{enumerate}
First, we easily show from the above proposition that 
%%%%%%%%%%%%  lemma %%%%%%%%%%%%%%
\begin{prop}\label{prop:6-3}
Let $f$ be a smooth function on $M$. Then 
\begin{equation}
[D_{+i},f]=p_{+i}((\mathrm{grad}f)^{1,0}),  \quad [D_{-i},f]=p_{-i}((\mathrm{grad}f)^{0,1}), \label{eqn:6-7}
\end{equation}
where $(\mathrm{grad}f)^{1,0}$ and $(\mathrm{grad}f)^{0,1}$ are the $(1,0)$- and $(0,1)$-part of gradient vector field of $f$.
\end{prop}
%%%%%%%
Next, we shall prove that the covariant derivative $\nabla$ on $\mathbf{S}_{\rho}$ is compatible with the ones on $T^{1,0}(M)$ and $T^{0,1}(M)$. We take a local unitary frame $\{\e_k\}_k$ and the connection $\omega$ on $\mathbf{U}(M)$. The covariant derivative $\nabla$ on $T^{1,0}(M)$ is expressed as 
\begin{equation}
\nabla_V\e_l=\sum_k \omega^k_l(V)\e_k=\sum_k g(\nabla_V\e_k,\bar{\e}_l)\e_k,
               \nonumber %\label{eqn:6-8}
\end{equation}
where $(\omega^k_l)_{k,l}$ is a local connection $1$-form of $\omega$ with respect to $\{\e_k\}_k$. Since the frame $\{\e_k\}_k$ induces a local unitary frame $\{s_{\alpha}\}_{\alpha=1}^{\dim V_{\rho}}$ of $\mathbf{S}_{\rho}$, the covariant derivative $\nabla$ on $\mathbf{S}_{\rho}$ is expressed as 
\begin{equation}
\nabla_V s_{\alpha}=\sum \omega^k_l(V) \pi_{\rho}(e_{kl})s_{\alpha}. 
  \nonumber%\label{eqn:6-9}
\end{equation}
Here, $e_{kl}$ is a local section of $\mathbf{U}(M)\times_{\Ad}\mathfrak{gl}(m,\mathbb{C})$ corresponding to a local section $\e_k\otimes \bar{\e}_l$ of $T^{1,0}(M)\otimes T^{0,1}(M)$. From this local expression of $\nabla$, we have
%%%%%%%%%%
\begin{prop}
The Clifford homomorphisms $p_{+i}(X)$ and $p_{-i}(\overline{X})$ for $X$ in $\Gamma (T^{1,0}(X))$ satisfy that 
\begin{equation}
\begin{split}
\nabla_V(p_{+i}(X)\phi)&=p_{+i}(\nabla_VX)\phi+p_{+i}(X)\nabla_V\phi, \\
\nabla_V(p_{-i}(\overline{X})\phi)&=p_{+i}(\nabla_V\overline{X})\phi+p_{+i}(\overline{X})\nabla_V\phi. \\
\end{split}
\end{equation}
for $\phi$ in $\Gamma (\mathbf{S}_{\rho})$ and $V$ in $\Gamma (T(M))$. 
\end{prop}
\begin{proof}
We take local frames $\{\e_k\}_k$ of $T^{1,0}(M)$ and $\{s_{\alpha}\}_{\alpha}$ of $\mathbf{S}_{\rho}$ as above.  Then 
\begin{equation}
\begin{split}
 &p_{+i}(\nabla_V \e_k)s_{\alpha}+p_{+i}(\e_k)\nabla_V s_{\alpha} \\ 
=&\sum_l \omega^l_k(V)p_{+i}(\e_l)s_{\alpha} +\sum \omega^s_t(V)p_{+i}(\e_k)\pi_{\rho}(e_{st})s_{\alpha}\\ 
=&\sum_l \omega^l_k(V)p_{+i}(\e_l)s_{\alpha} +\sum \omega^s_t(V)(-\delta_{kt}p_{+i}(\e_s)+\pi_{\rho+\mu_i}(e_{st})p_{+i}(\e_k))s_{\alpha} \\
=&\sum_l \omega^l_k(V)p_{+i}(\e_l)s_{\alpha} -\sum \omega^s_k(V)p_{+i}(\e_s) s_{\alpha}+\sum \omega^s_t(V) \pi_{\rho+\mu_i}(e_{st})p_{+i}(\e_k)s_{\alpha} \\
=&\nabla_V(p_{+i}(\e_k)s_{\alpha}),
\end{split}\nonumber
\end{equation}
where we use the infinitesimal expression of the equation \eqref{eqn:3-8},
$$
\delta_{kt}p_{+i}(\e_s)=\pi_{\rho+\mu_i}(e_{st})p_{+i}(\e_k)-p_{+i}(\e_k)\pi_{\rho}(e_{st}). 
$$
\end{proof}

In the rest of this section, we prove a conformal covariance of $D_{\pm i}$. The method is same to the conformal covariance of the Dirac operator (cf. \cite{LM}). The connection on $\mathbf{U}(M)$ in our definition of $D_{\pm i}$ may be arbitrary. But we should choose more geometric connections on $\mathbf{U}(M)$ as follows:
\begin{enumerate}
	\item  The Hermitian connection $\omega^h$ when $J$ is integrable.
	\item  The $\mathfrak{u}(m)$-component $\omega^u$ of the Levi-Civita connection $\omega^g$ with respect to $g$.
\end{enumerate}
The first connection $\omega^h$ is unique connection such that $\nabla^h g=0$, $\nabla^h J=0$, and the $(1,1)$-part of torsion tensor vanishs. We define the second connection $\omega^u$ more precisely. We take the Levi-Civita connection $\omega^g$, whose connection $1$-form is $\mathfrak{so}(2m)$-valued. Considering the inclusion $\mathfrak{u}(m)\subset \mathfrak{so}(2m)$, we can take the $\mathfrak{u}(m)$-component $\omega^u$ of $\omega^g$, which is a connection on $\mathbf{U}(M)$ (see \cite{KN}). Such a connection is given by
\begin{equation}
\nabla^u_V W=\nabla^g_V W-\frac{1}{2}J(\nabla^g_V J)(W)=\frac{1}{2}(\nabla^g_VW-J(\nabla^g_V(JW)))\label{eqn:6-10}
\end{equation}
for $V$ and $W$ in $\Gamma (T(M))$. 
\begin{proof}
For  $V$ in $T(M)$, the projections to $(1,0)$- and $(0,1)$-part are
\begin{equation}
  \pi^{1,0}(V)=\frac{1}{2}(1-\sqrt{-1}J)V \in T^{1,0}(M), \quad \pi^{0,1}(V)=\frac{1}{2}(1+\sqrt{-1}J)V \in T^{0,1}(M).\nonumber
\end{equation}
Then the $\mathfrak{u}(m)$-component $\omega^u$ of the Levi-Civita connection is\begin{equation}
\begin{split}
\nabla^u_VW &=\pi^{1,0}\nabla^g_V\pi^{1,0}W+\pi^{0,1}\nabla^g_V\pi^{0,1}W \\
      &=\frac{1}{2}(\nabla^g_VW-J(\nabla^g_V(JW))) \\
      &=\nabla^g_VW-\frac{1}{2}J(\nabla^g_VJ)(W).
\end{split}\nonumber
\end{equation}
\end{proof}
We can easily show that this covariant derivative $\nabla^u$ preserves $g$ and $J$ and gives a connection $\omega^u$ on $\mathbf{U}(M)$. But, the torsion tensor of $\omega^u$ does not always vanish. In fact, the torsion tensor is a $(1,1)$-tensor. We know that an almost Hermitian manifold $M$ is a K\"ahler manifold if and only if $\omega^u=\omega^h=\omega^g$. 

Let us consider a conformal deformation of Hermitian metric, $g\mapsto g'=e^{2\sigma (x)}g$ for $\sigma (x)$ in $C^{\infty}(M)$. The metric $g'$ is also a Hermitian metric on $M$. This deformation induces a principal bundle isomorphism 
\begin{equation}
\Psi:\mathbf{U}(M)\ni p=(\e_1,\cdots,\e_m)\mapsto p'=e^{-\sigma (x)}(\e_1,\cdots,\e_m)\in \mathbf{U}'(M), \nonumber
\end{equation}
and a bundle isomorphism 
\begin{equation}
\psi_{\rho}:\mathbf{S}_{\rho}=\mathbf{U}(M)\times_{\rho}V_{\rho}\ni[p,\phi]\mapsto [p',\phi]\in \mathbf{U}'(M)\times_{\rho}V_{\rho}=\mathbf{S}'_{\rho}. \nonumber
\end{equation}
Then the operators $D_{\pm i}$ and $D_{\pm i}'$ on $\mathbf{S}_{\rho}$ and $\mathbf{S}'_{\rho}$ are related as follows:
%%%%%%%%%%%%% theorem %%%%%%%%%%%%%%%%%%%%% 
\begin{thm}\label{thm:6-3}
\begin{enumerate}
	\item Let $M$ be a Hermitian manifold with Hermitian metric $g$. The K\"ahlerian gradients $\{D_{\pm i}\}_i$ associated to the Hermitian connection have conformal covariance, 
\begin{equation}
D'_{\pm i}=e^{(\pm \pi_{\rho}(c_1)-1) \sigma} \psi_{\rho\pm \mu_i} \circ D_{\pm i}\circ  (e^{\pm \pi_{\rho}(c_1)\sigma}\psi_{\rho})^{-1}.\label{eqn:6-11}
 \end{equation}
 Here, $c_1$ is the Casimir element with degree $1$. 
 %%%%%
 \item 	Let $M$ be an almost Hermitian manifold with Hermitian metric $g$. The operators $\{D_{\pm i}\}_i$ associated to the $\mathfrak{u}(m)$-component of the Levi-Civita connection have conformal covariance,
\begin{equation}
D'_{\pm i}=e^{(-w_{\pm i}-1) \sigma} \psi_{\rho\pm \mu_i} \circ D_{\pm i}\circ  (e^{-w_{\pm i}\sigma}\psi_{\rho})^{-1}. \label{eqn:6-12}
 \end{equation}
 Here, $w_{\pm i}$ is the conformal weight associated to $\rho$. 
\end{enumerate}
\end{thm}
%%%%%%%%%%%%%%%%%%%%%%%%%%%
\begin{rem}
\begin{enumerate}
	\item This theorem provides the reason why we call $w_{\pm i}$ the conformal weight.
	\item When $\dim \ker D_{\pm i}$ is finite, $\dim \ker D_{\pm i}$ is a conformal invariant of $M$.
\end{enumerate}
 
\end{rem}
\begin{proof}
First, we shall prove the conformal covariance \eqref{eqn:6-11}. We set $h_{k\bar{l}}:=g(\e_k,\bar{\e}_l)$ and $H:=(h_{k\bar{l}})_{kl}$, then a local connection $1$-from of $\omega^h$ is given by 
$$
H^{-1}\partial H=(h^{k\bar{s}}\partial h_{s\bar{l}})_{kl}.
$$
Under the conformal deformation $g\mapsto g'=e^{2\sigma}g$, we have 
\begin{equation}
\begin{split}
\omega'{}^k_l(V) &=g'(\nabla'_V\e'_k,\bar{\e}'_j)  \\
           &=g'(\nabla_V\e'_k+((V-\sqrt{-1}JV)\sigma)\e'_k,\e'_l) \\
           &=-\sqrt{-1}(JV\sigma)\delta_{kl}+\omega^k_l(V),
        \end{split}\nonumber
\end{equation}
and  
\begin{equation}
\begin{split}
\nabla'_V \psi_{\rho} (s_{\alpha})&=\sum \omega'{}^k_l (V) \pi_{\rho}(e'_{kl})\psi_{\rho} (s_{\alpha}) \\
 &=\psi_{\rho} \{ \sum (\omega^k_l(V)-\sqrt{-1}(JV\sigma)\delta_{kl})\pi_{\rho}(e_{kl})s_{\alpha}\} \\
 &=\psi_{\rho} \{ \nabla_V s_{\alpha}-\sqrt{-1}\pi_{\rho}(c_1) (JV\sigma)s_{\alpha} \}
\end{split}\nonumber
\end{equation}
on each associated bundle. Then, for $\phi$ in $\Gamma (\mathbf{S}_{\rho})$, 
\begin{equation}
\begin{split}
D'_{+i}\psi_{\rho} (\phi)&=\sum_k p_{+i}(\e'_k)\nabla'_{\bar{\e}'_k}\psi_{\rho} (\phi) \\
            &=e^{-\sigma}\psi_{\rho+\mu_i}\{ \sum_k p_{+i}(\e_k)(\nabla_{\bar{\e}_k}-\pi_{\rho}(c_1)(\bar{\e}_k\sigma))\phi\} \\
           &=e^{-\sigma}\psi_{\rho+\mu_i}\{ D_{+i}\phi-\pi_{\rho}(c_1)p_{+i}((\mathrm{grad}\sigma)^{1,0})\phi\}.
\end{split}\nonumber
\end{equation}
It follows from this equation and \eqref{eqn:6-7} that 
\begin{equation}
D'_{\pm i}=e^{(\pm \pi_{\rho}(c_1)-1) \sigma} \psi_{\rho\pm \mu_i} \circ D_{\pm i}\circ  (e^{\pm \pi_{\rho}(c_1)\sigma}\psi_{\rho})^{-1}.\nonumber
 \end{equation}
Next, we shall prove the conformal covariance \eqref{eqn:6-12}. The Levi-Civita connection changes under the conformal deformation $g\mapsto g'=e^{2\sigma}g$ as $$
\nabla^{g'}_VW=\nabla^g_VW+(V\sigma)W+(W\sigma)V-g(V,W)\sigma.
$$
So the local connection $1$-from of the $\mathfrak{u}(m)$-component changes as
\begin{equation}
\omega'{}^k_l(V)=\omega^k_l(V)+(\e_l\sigma)g(V,\bar{\e}_k)-(\bar{\e}_k\sigma)g(V,\e_l),\nonumber
\end{equation}
 and the covariant derivative on each associated bundle does as 
$$
\nabla'_{V}\psi_{\rho} (s_{\alpha}) =\psi_{\rho} \{ \nabla_V s_{\alpha}+\sum (\e_l\sigma)g(V,\bar{\e}_k)\pi_{\rho}(e_{kl})s_{\alpha}-\sum (\bar{\e}_k\sigma)g(V,\e_l)\pi_{\rho}(e_{kl}) s_{\alpha}\}.
$$
Then, from \eqref{eqn:3-16}, we have
\begin{equation}
\begin{split}
D'_{+i}\psi_{\rho} (\phi) &=e^{-\sigma}\psi_{\rho+\mu_i} \{\sum_k p_{+i}(\e_k)(\nabla_{\bar{\e}_k}-\sum_{i}(\bar{\e}_l\sigma)\pi_{\rho}(e_{lk}))\phi \} \\
&=e^{-\sigma} \psi_{\rho+\mu_i} \{D_{+i}\phi +w_{+i}p_{+i}((\mathrm{grad}\sigma)^{1,0}) \phi \}    
\end{split}\nonumber
\end{equation}
Then we have the conformal covariance \eqref{eqn:6-12}. Similarly we can prove the conformal invariance for $D_{-i}$. 
\end{proof}
%%%%%%%%%%%%%%%  rem %%%%%%%%%%%%%%%%
\begin{exam}
We consider a Riemannian surface $M$ whose Riemannian metric is always a K\"ahler metric. The irreducible representations of $U(1)$ are parametrized by $l$ in $\mathbb{Z}$. For each irreducible representation $(\pi_{l},V_l)$ , we have the conformal weight $w_{\pm 1}=\mp l$. On the other hands, we know that the Casimir operator $\pi_{l}(c_1)$ is $l$ on $V_l$. Thus, the conformal covariance \eqref{eqn:6-11} coincides with \eqref{eqn:6-12} on Riemannian surfaces. 
\end{exam}
%%%%%%%%%%
%%%%%%%%%%%%%%%%%%%%%%%%%%%%%%%%%%%%%%%%%%%%%%%%%%
%%%%%%%%%%%%%%%%%%%%%%%%%%%%%%%%
%%           7               %%
%%%%%%%%%%%%%%%%%%%%%%%%%%%%%%%%
\section{The Bochner identities on K\"ahler manifolds}\label{sec:7}
In this section, we extend the relations for the Clifford homomorphisms to the (generalized) Bochner identities for the K\"ahlerian gradients on K\"ahler manifolds. From now on, we assume that $M$ is a K\"ahler manifold with K\"ahler metric $g$, and the connection on $\mathbf{U}(M)$ is the Levi-Civita connection. 

We define a second order differential operator $\nabla_{V,W}^2$ on $\mathbf{S}_{\rho}$ for vector fields $V$ and $W$ by 
\begin{equation}
\nabla_{V,W}^2:=\nabla_V\nabla_W-\nabla_{\nabla_VW}:\Gamma (\mathbf{S}_{\rho})
\to \Gamma (\mathbf{S}_{\rho}). \label{eqn:7-1}
\end{equation}
We denote the formal adjoint operator of $\nabla$ by $\nabla^{\ast}$. Then the connection Laplacian is given by 
\begin{equation}
\nabla^{\ast}\nabla=-\sum_k (\nabla_{e_k,e_k}^2+\nabla_{Je_k,Je_k}^2)=-\sum \nabla_{\bar{\e}_k,\e_k}^2-\sum \nabla_{\e_k,\bar{\e}_k}^2. \label{eqn:7-2}
\end{equation}
Thus the connection Laplacian splits into the sum of $\nabla^{1,0}{}^{\ast}\nabla^{1,0}$ and $\nabla^{0,1}{}^{\ast}\nabla^{0,1}$,
\begin{equation}
\nabla^{1,0}{}^{\ast}\nabla^{1,0}=-\sum \nabla^2_{\bar{\e}_k,\e_k}, \quad \nabla^{0,1}{}^{\ast}\nabla^{0,1}=-\sum \nabla^2_{\e_k,\bar{\e}_k}. \label{eqn:7-3}
\end{equation}

On the other hand, the difference between $\nabla^{1,0}{}^{\ast}\nabla^{1,0}$ and $\nabla^{0,1}{}^{\ast}\nabla^{0,1}$ gives a bundle endomorphism of $\mathbf{S}_{\rho}$ depending on the curvature. Furthermore, we show that some linear combinations of $\{D_{\pm i}\}_i$ are also bundle endomorphisms on $\mathbf{S}_{\rho}$ depending on the curvature. So, for a while, we discuss such endomorphisms.

The curvature $R_T$ on tangent bundle of $T(M)$ is given by
\begin{equation}
R_T(V,W):=\nabla^2_{V,W}-\nabla^2_{W,V} \quad \textrm{for $V,W$ in $T(M)$}. \label{eqn:7-4}
\end{equation}
On K\"ahler manifolds, this curvature tensor $R_T$ satisfies $J R_T(V,W)=R_T(V,W) J$ and $R_T(JV,JW)=R_T(V,W)$. So, $R_T$ is an $\mathfrak{u}(m)$-valued $(1,1)$-form and leads to the curvature on $T^{1,0}(M)$ and $T^{0,1}(M)$. For example, the curvature on $T^{1,0}(M)$ is 
\begin{equation}
\begin{split}
R_{T^{1,0}}(V,W) &=\sum g(R_T(V,W)\e_l,\bar{\e}_k)\e_k\otimes \bar{\e}_l \\
                 &=\sum g(R_T(V,W)\e_l,\bar{\e}_k)e_{kl},
 \end{split} \label{eqn:7-5}
\end{equation}
where $g$ is extended complex linearly to a complex metric on $T(M)\otimes \mathbb{C}$. 

Contracting components of $R_T$, we have the Ricci curvature $\mathrm{Ric}$ and the scalar curvature $\kappa$, 
\begin{gather}
\mathrm{Ric}(V,W)=\sqrt{-1} g(\sum R_T(\e_k,\bar{\e}_k)V,JW), 
                                           \label{eqn:7-6}   \\
\kappa=2\sum g(R_T(\e_k,\bar{\e}_k)\e_l,\bar{\e}_l). \label{eqn:7-7}
\end{gather}

Let $R_{\rho}$ be the curvature of $\nabla$ on $\mathbf{S}_{\rho}$,
\begin{equation}
R_{\rho}(V,W):=\nabla_V\nabla_W-\nabla_W\nabla_V-\nabla_{[V,W]}=\nabla^2_{V,W}-\nabla^2_{W,V}. \label{eqn:7-8}
\end{equation}
Since the derivative $\nabla$ on $\mathbf{S}_{\rho}$ is induced from the Levi-Civita connection, we have a formula of $R_{\rho}$,
\begin{equation}
R_{\rho}(V,W)=\sum g(R_T(V,W)\e_l,\bar{\e}_k)\pi_{\rho}(e_{kl}). \label{eqn:7-9}
\end{equation}

We shall contract the curvature $R_{\rho}$ by the action of the enveloping algebra bundle on $M$ and give some bundle endomorphisms of $\mathbf{S}_{\rho}$. We define the enveloping algebra bundle $\mathbf{U}(\gl)$ on $M$ to be
\begin{equation}
\mathbf{U}(\gl):=\mathbf{U}(M)\times_{\Ad} U(\gl). \label{eqn:7-10}
\end{equation}
Here, $\Ad$ is the adjoint representation of $U(m)$ on $U(\gl)$. Each fiber of $\mathbf{U}(\gl)$ has an algebra structure isomorphic to $U(\gl)$, and each associated bundle $\mathbf{S}_{\rho}$ is a bundle of $\mathbf{U}(\gl)$-module. 
%%%%%%%%%%%%%%% definition %%%%%%%%%%%%%%%%%%%
\begin{defini}\label{def:7-1}
Let $R_{\rho}$ be the curvature of $\nabla$ on $\mathbf{S}_{\rho}$. We define a set of bundle endomorphisms $\{R_{\rho}^q,\tilde{R}_{\rho}^q \}_{q\ge 0}$ depending on $R_{\rho}$ by
\begin{gather}
R_{\rho}^q:=\sum_{kl} \pi_{\rho}(e_{lk}^q)R_{\rho}(\e_k,\bar{\e}_l) \in \Gamma (\End (\mathbf{S}_{\rho})), \label{eqn:7-11}\\
\tilde{R}_{\rho}^q:=\sum_{kl} \pi_{\rho}(\tilde{e}_{kl}^q)R_{\rho}(\e_k,\bar{\e}_l) \in \Gamma (\End (\mathbf{S}_{\rho})),\label{eqn:7-12}
\end{gather}
where $e_{kl}^q$ and $\tilde{e}_{kl}^q$ are local sections of $\mathbf{U}(\gl)$ corresponding to \eqref{eqn:2-1} and \eqref{eqn:2-9}. We call these bundle endomorphisms \textit{the curvature endomorphisms} on $\mathbf{S}_{\rho}$. 
\end{defini}
%%%%%%%%%%%%%%%

The curvature endomorphisms have real eigenvalues on each point. In fact, we have
%%%%%%%%%%%%%%%%%%% prop %%%%%%%%%%%%%%%%%%%%%
\begin{prop}\label{prop:7-2}
The curvature endomorphisms $R_{\rho}^q$ and $\tilde{R}_{\rho}^q$ are self adjoint bundle endomorphisms, that is, 
\begin{equation}
(R_{\rho}^q\phi_1,\phi_2)_x=(\phi_1,R_{\rho}^q \phi_2)_x, \quad (\tilde{R}_{\rho}^q\phi_1,\phi_2)_x=(\phi_1,\tilde{R}_{\rho}^q \phi_2)_x  
\end{equation}
for $x$ in $M$ and $\phi_1,\phi_2$ in $(\mathbf{S}_{\rho})_x$. Furthermore, it holds that
\begin{equation}
\sum_{p=0}^q \binom{q}{p}(-m)^{q-p}\tilde{R}_{\rho}^p=(-1)^q\sum_{p=0}^q \pi_{\rho}(K_{q-p}(-c))R_{\rho}^q. \label{eqn:7-15}
\end{equation}
The curvature endomorphism $R_{\rho}^q$ is called positive (resp. negative) if all the eigenvalues of $(R_{\rho}^q)_x$ is positive (resp. negative) for any $x$ in $M$. 
\end{prop}
%%%%%%%%%%%%%%
%%%%%%%%%% proof %%%%%%%%%%%%%%
\begin{proof}
From the definition of $R_{\rho}^q$, we have 
$$
R^q_{\rho}=\sum g(R_T(\e_i,\bar{\e}_j)\e_l,\bar{\e}_k)   \pi_{\rho}(e_{ji}^q)\pi_{\rho}(e_{kl})=\sum g(R_T(\e_i,\bar{\e}_j)\e_l,\bar{\e}_k)\pi_{\rho}(e_{ji}^qe_{kl}).
$$
Then 
\begin{equation}
\begin{split}
 &R_{\rho}^q \\
= &\sum g(R_T(\e_i,\bar{\e}_j)\e_l,\bar{\e}_k) \pi_{\rho}(e_{ji}^qe_{kl}) \\
= &\sum g(R_T(\e_i,\bar{\e}_j)\e_l,\bar{\e}_k)  (\pi_{\rho}(e_{kl}e_{ji}^q)-\delta_{lj}\pi_{\rho}(e_{ki}^q)+\delta_{ki}\pi_{\rho}(e_{jl}^q))\quad \textrm{(by \eqref{eqn:2-2})} \\
= &\sum g(R_T(\e_i,\bar{\e}_j)\e_l,\bar{\e}_k) \pi_{\rho}(e_{kl}e_{ji}^q)-g(R_T(\e_i,\bar{\e}_j)\e_j,\bar{\e}_k)\pi_{\rho}(e_{ki}^q)+ g(R_T(\e_i,\bar{\e}_j)\e_l,\bar{\e}_i)\pi_{\rho}(e_{jl}^q) \\
= &\sum g(R_T(\e_i,\bar{\e}_j)\e_l,\bar{\e}_k) \pi_{\rho}(e_{kl}e_{ji}^q)- g(R_T(\e_l,\bar{\e}_l)\e_i,\bar{\e}_k)\pi_{\rho}(e_{ki}^q)+ g(R_T(\e_i,\bar{\e}_i)\e_l,\bar{\e}_j)\pi_{\rho}(e_{jl}^q) \\
=&\sum g(R_T(\e_i,\bar{\e}_j)\e_l,\bar{\e}_k) \pi_{\rho}(e_{kl}) \pi_{\rho}(e_{ji}^q),
\end{split}\nonumber
\end{equation}
where we use the Bianchi's identities for $R_T$. Then we have
\begin{equation}
\begin{split}
 &(R^q_{\rho})^{\ast}= \sum g(R(\bar{\e}_i,\e_j)\bar{\e}_l,\e_k) \pi_{\rho}(e_{kl})^{\ast}\pi_{\rho}(e_{i_{q-1}i})^{\ast}\cdots \pi_{\rho}(e_{ji_1})^{\ast} \\
=& \sum g(R(\e_j,\bar{\e}_i)\e_k,\bar{\e}_l) \pi_{\rho}(e_{lk}) \pi_{\rho}(e_{ii_{q-1}})\cdots \pi_{\rho}(e_{i_1j}) \\
=&\sum g(R(\e_i,\bar{\e}_j)\e_l,\bar{\e}_k) \pi_{\rho}(e_{kl}) \pi_{\rho}(e_{ji_{q-1}})\cdots \pi_{\rho}(e_{i_1i}) =R^q_{\rho}.
\end{split}\nonumber
\end{equation}
The second assertion can be proved by using theorem \ref{thm:4-3}. 
\end{proof}
%%%%%
%%%%%%%%%%% example %%%%%%%%%%%%%%%%%%
\begin{exam}\label{ex:7-1}
In the case of $q=0$, the curvature endomorphisms $R_{\rho}^0$ and $\tilde{R}_{\rho}^0$ become
\begin{equation}
R_{\rho}^0=\tilde{R}_{\rho}^0=\sum R_{\rho}(\e_k,\bar{\e}_k)=\sum_{i,j} g(\sum_k R_{T}(\e_k,\bar{\e}_k)\e_j,\bar{\e}_i)\pi_{\rho}(e_{ij}).\label{eqn:7-13}
\end{equation}
This is the so-called mean curvature on $\mathbf{S}_{\rho}$. We show that, if $M$ is a Ricci-flat K\"ahler manifold, then $R_{\rho}^0=\tilde{R}_{\rho}^0=0$. 
\end{exam}
%%%%%%%
%%%%%%%%%% example %%%%%%%%%%%%%%%%%%
\begin{exam}\label{ex:7-2}
If $M$ is a K\"ahler manifold of constant holomorphic sectional curvature$=r$, then $R_T$ on $T(M)$ satisfies 
$$
g(R_T(\e_i,\bar{\e}_j)\e_l,\bar{\e}_k)=\frac{r}{2}(\delta_{kl}\delta_{ij}+\delta_{lj}\delta_{ki}).
$$
The curvature endomorphism $R^q_{\rho}$ is also constant for any $q$, which is given by
\begin{equation}
R^q_{\rho}=\frac{r}{2}\pi_{\rho}(c_qc_1+c_{q+1}). \label{eqn:7-14}
\end{equation}
\end{exam}
%%%%%%%
%%%%%%% example %%%%%%%%%%%%%%%%%
\begin{exam}\label{ex:7-3}
We think of the curvature endomorphisms on the spinor bundle $\sum \Lambda^{0,p}(M)$, where $\Lambda^{0,p}(M)$ is $\mathbf{S}_{\Lambda^{0,p}}$ with $\Lambda^{0,p}:=(1_p,0_{m-p})$. It is from \eqref{eqn:5-5} that the $0$-th curvature endomorphism is 
\begin{equation}
R_{\Lambda^{0,p}}^0=\sum g(R_{T}(\e_k,\bar{\e}_k)\e_j,\bar{\e}_i)\pi_{\Lambda^{0,p}}(e_{ij})=-\sum g(R_{T}(\e_k,\bar{\e}_k)\e_j,\bar{\e}_i) \e_i\cdot\bar{\e}_j\cdot.\nonumber
\end{equation}
In particular, we have
\begin{equation}
R_{\Lambda^{0,m}}=\sum g(R_{T}(\e_k,\bar{\e}_k)\e_j,\bar{\e}_i)\delta_{ij}=\frac{\kappa}{2}.\nonumber
\end{equation}
We calculate the first curvature endomorphism $R_{\Lambda^{0,p}}^1$. From the definition of $R_{\Lambda^{0,p}}^1$, we have
\begin{equation}
\begin{split}
R_{\Lambda^{0,p}}^1&=\sum g(R_{T}(\e_k,\bar{\e}_l)\e_j,\bar{\e}_i)\pi_{\Lambda^{0,p}}(e_{lk})\pi_{\Lambda^{0,p}}(e_{ij}) \\
   &=\sum g(R_{T}(\e_k,\bar{\e}_l)\e_j,\bar{\e}_i)\e_l\cdot\bar{\e}_k \cdot\e_i\cdot\bar{\e}_j\cdot.   \end{split}\nonumber
\end{equation}
Then the Bianchi's identity and the Clifford algebra relations lead us to $R^1_{\Lambda^{0,p}}=R^0_{\Lambda^{0,p}}$. Thus, we conclude that 
\begin{equation}
R^1_{\Lambda^{0,p}}=R_{\Lambda^{0,p}}^0=
\begin{cases}
0 &  \textrm{for $p=0$}, \\
\kappa/2 &\textrm{for $p=m$}, \\
\sum -g(R_{T}(\e_k,\bar{\e}_k)\e_j,\bar{\e}_i)\e_i\cdot\bar{\e}_j\cdot & \textrm{otherwise}.\end{cases}\label{eqn:cur}
\end{equation}
The curvature endomorphisms $R^1_{\Lambda^{0,p}}$ and $R^0_{\Lambda^{0,p}}$ depend on the Ricci curvature of $M$. The Ricci curvature is called positive if $\mathrm{Ric}(V,V)>0$ for any real vector field $V\neq 0$, or $\mathrm{Ric}(X,\bar{X})>0$ for any $(1,0)$-vector field $X\neq 0$. Then we show that, if the Ricci tensor is positive (resp. negative), then $R_{\Lambda^{0,p}}^1=R^0_{\Lambda^{0,p}}$ is positive (resp. negative) except $p=0$. In fact, we take eigenbasis $\{(\e_k)_x\}_{k=1}^m$ of $\mathrm{Ric}_x$ such that 
$$\mathrm{Ric}((\e_k)_x,(\bar{\e}_l)_x)=\mathrm{Ric}((e_k)_x,(e_l)_x)=\mathrm{Ric}(J_x(e_k)_x,J_x(e_l)_x)=\lambda_k\delta_{kl}.$$
Then the eigenvalues of $(R_{\Lambda^{0,p}}^0)_x$ are $\{\lambda_{i_1}+\cdots+\lambda_{i_p}| 1\le i_1<\cdots <i_p\le m\}$. 
\end{exam}
%%%%%%%%%%%

We consider the formal adjoint operator of $D_{\pm i}$, whose principal symbol is the Clifford homomorphism $-p_{\pm i}^{\ast}$. 
%%%%%%%%%%%%%%%  proposition %%%%%%%%%%%%%%%%
\begin{prop}\label{prop:7-3}
We denote by $D_{\pm i}^{\ast}$ the formal adjoint operator of $D_{\pm i}$ such that 
\begin{equation}
\int_M (D_{\pm i}\phi,\psi)dx=\int_M (\phi,D_{\pm i}^{\ast}\psi)dx
\end{equation}
for all compactly supported $\phi$ in $\Gamma (\mathbf{S}_{\rho})$ and $\psi$ in $\Gamma (\mathbf{S}_{\rho\pm i})$. 
Then $D_{\pm i}^{\ast}$ is expressed as follows. 
\begin{gather}
D_{+i}^{\ast}=-\sum p_{+i}(\e_k)^{\ast}\nabla_{\e_k}:\Gamma (\mathbf{S}_{\rho+\mu_i})\to \Gamma (\mathbf{S}_{\rho}),\label{eqn:7-16}\\
D_{-i}^{\ast}=-\sum p_{-i}(\bar{\e}_k)^{\ast}\nabla_{\bar{\e}_k}:\Gamma (\mathbf{S}_{\rho-\mu_i})\to \Gamma (\mathbf{S}_{\rho}).\label{eqn:7-17}
\end{gather}
\end{prop}
%%%%%%%%%
%%%%%%%%%%%% proof %%%%%%%%%%%%%
\begin{proof}
We fix $x$ in $M$ and choose a local Hermitian frame $\{\e_k\}_k$ in a neighborhood of $x$ such that $(\nabla\e_k)_x=0$ for any $k$. Then, at $x$, 
\begin{equation}
\begin{split}
(D_{+i}\phi,\psi)=&\sum (p_{+i}(\e_k)\nabla_{\bar{\e}_k}\phi,\psi)=\sum (\nabla_{\bar{\e}_k}\phi,p_{+i}(\e_k)^{\ast}\psi)\\
                 =&\sum \{\bar{\e}_k(\phi,p_{+i}(\e_k)^{\ast}\psi)-(\phi,p_{+i}(\nabla_{\bar{\e}_k}\e_k)^{\ast} \psi)-(\phi,p_{+i}(\e_k)^{\ast} \nabla_{\e_k}\psi)\} \\
=&\sum \bar{\e}_k(\phi,p_{+i}(\e_k)^{\ast}\psi)    +(\phi,-\sum p_{+i}(\e_k)^{\ast}\nabla_{\e_k}\psi) \\
=& \mathrm{div} (\bar{X}) +(\phi,-\sum p_{+i}(\e_k)^{\ast}\nabla_{\e_k}\psi),
\end{split}\nonumber
\end{equation}
where $\bar{X}$ is $(0,1)$-vector field defined by the condition that $g(\bar{X},Y)=( \phi,p_{+i}(Y)^{\ast}\psi)$ for any $(1,0)$-vector $Y$. We show that, on K\"ahler manifolds, the divergence of $(0,1)$-vector field $\bar{X}$ (resp. $(1,0)$-vector field $X$) is $\mathrm{div}(\bar{X})=\sum g(\nabla_{\bar{\e}_k}\bar{X}, \e_k)$ (resp $\mathrm{div}(X)=\sum g(\nabla_{\e_k}X, \bar{\e}_k)$). So we have proved the proposition. 
\end{proof}
%%%%%%%%
\begin{lem}
The second order differential operators $\{D_{\pm i}^{\ast}D_{\pm i}\}_i$ satisfy that 
\begin{gather}
\sum_{1\le i\le m} w_{+i}^q D_{+i}^{\ast}D_{+i}=-\sum \pi_{\rho}(\tilde{e}_{kl}^q)\nabla^2_{\e_k,\bar{\e}_l}, \label{eqn:7-18}\\
\sum_{1\le i\le m} w_{-i}^q D_{-i}^{\ast}D_{-i}=-\sum \pi_{\rho}(e_{kl}^q)\nabla^2_{\bar{\e}_k,\e_l}. \label{eqn:7-19}
\end{gather}
\end{lem}
%%%%%%%%
\begin{proof}
This lemma is easily shown from \eqref{eqn:3-18} and \eqref{eqn:3-19}: 
\begin{equation}
\begin{split}
 \sum w_{+i}^q D_{+i}^{\ast}D_{+i}=&-\sum w_{+i}^q p_{+i}(\e_k)^{\ast}\nabla_{\e_k} p_{+i}(\e_l)\nabla_{\bar{\e}_l} \\
 =&-\sum w_{+i}^q p_{+i}(\e_k)^{\ast}p_{+i}(\e_l)\nabla_{\e_k}\nabla_{\bar{\e}_l}\\
 =&-\sum w_{+i}^q p_{+i}(\e_k)^{\ast}p_{+i}(\e_l)(\nabla_{\e_k}\nabla_{\bar{\e}_l}-\nabla_{\nabla_{\e_k}\bar{\e}_l})\\ 
=&-\sum w_{+i}^q p_{+i}(\e_k)^{\ast}p_{+i}(\e_l)\nabla_{\e_k,\bar{\e}_l}^2 =-\sum \pi_{\rho}(\tilde{e}_{kl}^q) \nabla^2_{\e_k,\bar{\e}_l}.
\end{split}\nonumber
\end{equation}
\end{proof}

We are ready for the Bochner identities for the K\"ahlerian gradients. 
%%%%%%%%% theorem %%%%%%%%%%%%%%%%%%%
\begin{thm}[The Bochner identities]\label{thm:7-4}
We consider the K\"ahlerian gradients $D_{\pm i}:\Gamma (\mathbf{S}_{\rho})\to \Gamma (\mathbf{S}_{\rho\pm \mu_i})$ and their formal adjoint operators $D_{\pm i}^{\ast}:\Gamma (\mathbf{S}_{\rho\pm \mu_i})\to \Gamma (\mathbf{S}_{\rho})$ on a K\"ahler manifold $M$. Then the operators satisfy the following identities:
\begin{enumerate}
	\item (in the case of degree $0$)
	\begin{gather}
\sum_{1\le i\le m}D_{-i}^{\ast}D_{-i}=\nabla^{1,0}{}^{\ast}\nabla^{1,0},	\quad \sum_{1\le i\le m}D_{+i}^{\ast}D_{+i}=\nabla^{0,1}{}^{\ast}\nabla^{0,1}, \label{eqn:7-20} \\
	\sum_{1\le i\le m}D_{-i}^{\ast}D_{-i}+D_{+i}^{\ast}D_{+i}=\nabla^{1,0}{}^{\ast}\nabla^{1,0}+\nabla^{0,1}{}^{\ast}\nabla^{0,1}=\nabla^{\ast}\nabla, \label{eqn:7-21} \\
	\sum_{1\le i\le m}D_{-i}^{\ast}D_{-i}-D_{+i}^{\ast}D_{+i}=\nabla^{1,0}{}^{\ast}\nabla^{1,0}-\nabla^{0,1}{}^{\ast}\nabla^{0,1}=R_{\rho}^0.\label{eqn:7-22} 
	\end{gather}
	In particular, if $M$ is a Ricci-flat K\"ahler manifold, then
	\begin{equation}
\sum_{1\le i\le m}D_{-i}^{\ast}D_{-i}=\sum_{1\le i\le m} D_{+i}^{\ast}D_{+i}.\label{eqn:7-23} 
	\end{equation}
	\item (in the case of degree $1$)
	\begin{equation}
	\sum_{1\le i\le m}w_{-i}D_{-i}^{\ast}D_{-i}+w_{+i}D_{+i}^{\ast}D_{+i}=R^1_{\rho}. \label{eqn:7-24} 
	\end{equation}
	\item (in the case of degree $q$)
\begin{multline}
 \sum_{0\le p\le q} \bp q \\ p\ep (-m)^{q-p}R^p_{\rho} 
=(-1)^{q+1} \sum_{1\le i\le m} (\sum_{0\le p\le q} \pi_{\rho}(K_{q-p}(-\tilde{c}))w_{+i}^p)D_{+i}^{\ast}D_{+i}\\+ \sum_{1\le i\le m} (w_{-i}-m)^q D_{-i}^{\ast}D_{-i}.\label{eqn:7-25} 
\end{multline}
\end{enumerate}
In particular, $D_{-i}^{\ast}D_{-i}$ is expressed as a linear combination of $\{D_{+j}^{\ast}D_{+j}\}_{j=1}^m$ and $\{ R^q_{\rho}\}_{q=0}^{m-1}$.
\end{thm}
%%%%%%%%%%
%%%%%%% remark %%%%%%%%%%
\begin{rem}
It follows from the discussion in remark \ref{rem:4-1} that, if there are $2N$ K\"ahlerian gradients on $\mathbf{S}_{\rho}$, then at most $N$ identities in \eqref{eqn:7-25} make sense. The identities with degree $q\ge N$ reduces to the ones for $0\le q \le N-1$. We can also have 
\begin{equation}
\begin{split}
&\sum_{0\le p\le q} \bp q \\ p\ep (-m)^{q-p}\tilde{R}^p_{\rho} \\
=&(-1)^{q+1} \sum_{1\le i\le m } (\sum_{0\le p\le q} \pi_{\rho}(K_{q-p}(-c))w_{-i}^p)D_{-i}^{\ast}D_{-i}+ \sum_{1\le i\le m } (w_{+i}-m)^q D_{+i}^{\ast}D_{+i} \end{split}\nonumber
\end{equation}
for any $q$. But, these identities reduce to \eqref{eqn:7-25}. 
\end{rem}
%%%%%%%%
%%%%%%%%%%% proof %%%%%%%%%%%%%%%%
\begin{proof}
We have known the relations \eqref{eqn:4-14} for principal symbols of the K\"ahlerian gradients, namely, the Clifford homomorphisms. Then 
\begin{equation}
\begin{split}
 &\sum_p \bp q \\ p\ep (-m)^{q-p}R^p_{\rho} \\
=&\sum_{p,k,l} \bp q \\ p\ep (-m)^{q-p} \pi_{\rho}(e_{lk}^p)(\nabla^2_{\e_k,\bar{\e}_l}-\nabla^2_{\bar{\e}_l,\e_k}) \\ 
=&\sum_{k,l,i} (w_{-i}-m)^q p_{-i}(\bar{\e}_l)^{\ast}p_{-i}(\bar{\e}_k)\nabla^2_{\e_k,\bar{\e}_l} +\sum_{i} (w_{-i}-m)^q D_{-i}^{\ast}D_{-i} \\
=&\sum_{k,l,i,p} (-1)^q \pi_{\rho}(K_{q-p}(-\tilde{c}))w_{+i}^p p_{+i}(\e_k)^{\ast}p_{+i}(\e_l)\nabla^2_{\e_k,\bar{\e}_l} +\sum_i (w_{-i}-m)^q D_{-i}^{\ast}D_{-i}\\
=&\sum_i (w_{-i}-m)^q D_{-i}^{\ast}D_{-i}+(-1)^{q+1} \sum_{i,p} \pi_{\rho}(K_{q-p}(-\tilde{c}))w_{+i}^pD_{+i}^{\ast}D_{+i}.
\end{split}\nonumber
\end{equation}
\end{proof}
%%%%%%%%%%%%%%

We sometimes deal with a vector bundle $\mathbf{S}_{\rho}\otimes E$, where $E$ is a Hermitian holomorphic vector bundle on $M$ with Hermitian connection $\nabla^E$. The Clifford homomorphisms and the K\"ahlerian gradients can be defined on $\mathbf{S}_{\rho}\otimes E$ like the twisted Dirac operator. The Clifford homomorphisms on $\mathbf{S}_{\rho}\otimes E$ are given by
$$
p_{+i}(X)\otimes \id :\mathbf{S}_{\rho}\otimes E\to \mathbf{S}_{\rho+\mu_i}\otimes E, \quad p_{+i}(\bar{X})\otimes \id :\mathbf{S}_{\rho}\otimes E\to \mathbf{S}_{\rho-\mu_i}\otimes E,
$$
and the operators are 
\begin{gather}
D_{+i}^E:=\sum_k (p_{+i}(\e_k)\otimes \id)(\nabla_{\bar{\e}_k}\otimes \id +\id \otimes \nabla^E_{\bar{\e}_k}), \\
D_{-i}^E:=\sum_k (p_{-i}(\bar{\e}_k)\otimes \id)(\nabla_{\e_k}\otimes \id +\id \otimes \nabla^E_{\e_k}) .
\end{gather}
We call these operators \textit{the twisted K\"ahlerian gradients} on twisted associated bundle. To provide the Bochner identities, we consider curvature endomorphisms on $\mathbf{S}_{\rho}\otimes E$ depending on the curvature for $E$. We take the curvature $R_E$ of $\nabla^E$ and define curvature endomorphisms $\mathfrak{R}_E^q$ by
\begin{equation}
\mathfrak{R}_E^q:=\sum \pi_{\rho}(e_{lk}^q)\otimes R_E(\e_k,\bar{\e}_l) \in \Gamma (\End (\mathbf{S}_{\rho}\otimes E)).
\end{equation}
In particular, $\mathfrak{R}_E^0$ depends on the mean curvature of $\nabla^E$, that is,
\begin{equation}
\mathfrak{R}^0_E=\id \otimes \sum R_{E}(\e_k,\bar{\e}_k).
\end{equation}
Then we have the Bochner identities for the twisted K\"ahlerian gradients.
%%%%%%%%%%%%%% corollary %%%%%%%%%%%%%%%%%%%
\begin{cor}
The twisted K\"ahlerian gradients $D^E_{\pm i}$ and $(D_{\pm i}^E)^{\ast}$ satisfy 
\begin{multline}
 \sum_{0\le p\le q} \bp q \\ p\ep (-m)^{q-p}(R^p_{\rho}+\mathfrak{R}_E^p) \\
=(-1)^{q+1} \sum_{1\le i\le m} (\sum_{0\le p\le q} \pi_{\rho}(K_{q-p}(-\tilde{c}))w_{+i}^p)(D_{+i}^E)^{\ast}D_{+i}^E+ \sum_{1\le i\le m} (w_{-i}-m)^q (D_{-i}^E)^{\ast}D_{-i}^E\label{eqn:7-26} 
\end{multline}
for any $q$ in $\mathbb{Z}_{\ge 0}$. Here, $R^q_{\rho}:=R^q_{\rho}\otimes \id$ depends on the curvature for $\mathbf{S}_{\rho}$, and $\mathfrak{R}_E^q$ does on the curvature for $E$.  
\end{cor}
%%%%%%%%%%%%%%%%
%%%%%%%%%%%%%%%%%%%%%%%%%%%%%%%%
%%           8               %%
%%%%%%%%%%%%%%%%%%%%%%%%%%%%%%%%
\section{Some applications}\label{sec:8}
Some applications are presented in this section, where we assume that $M$ is a closed K\"ahler manifold.

%%%%%%%%%%%% subsection %%%%%%%%%%%%%%%%
\subsection{Holomorphic sections and anti-holomorphic sections}
We discuss holomorphic structure on $\mathbf{S}_{\rho}$. Let $\mathbf{GL}_{\mathbb{C}}(M)$ be the principal bundle of holomorphic frames of K\"ahler manifold $M$ whose structure group is $GL(m,\mathbb{C})$. Since there is one-to-one correspondence between the complex representations of $GL(m,\mathbb{C})$ and the unitary representations of $U(m)$, the vector bundle $\mathbf{S}_{\rho}$ associated to $\mathbf{U}(M)$ is also associated to $\mathbf{GL}_{\mathbb{C}}(M)$, 
$$
\mathbf{S}_{\rho}=\mathbf{U}(M)\times_{\rho}V_{\rho}=\mathbf{GL}_{\mathbb{C}}(M)\times_{\rho}V_{\rho}.
$$
Thus, we know that the vector bundle $\mathbf{S}_{\rho}$ has a holomorphic Hermitian structure. A point to notice in holomorphic category is that we can not use a Hermitian frame, but a holomorphic frame. So the covariant derivative is expressed as 
\begin{equation}
\nabla^{1,0}\phi=\sum \nabla_{\partial/\partial z_i}\phi\otimes dz^i, \quad \nabla^{0,1}\phi=\sum \nabla_{\partial/\partial\bar{z}_i}\phi\otimes d\bar{z}^i,
 \nonumber
\end{equation}
and $D_{\pm i}$ is 
\begin{equation}
D_{-i}=\sum p_{-i}(dz^i)\nabla_{\partial/\partial z_i},\quad D_{+i}=\sum p_{+i}(d\bar{z}_i)\nabla_{\partial/\partial\bar{z}_i}. \nonumber
\end{equation}
The holomorphic and anti-holomorphic sections of $\mathbf{S}_{\rho}$ are characterized as 
\begin{equation}
\begin{split}
\ker \nabla^{0,1}&=\textrm{\{holomorphic sections of $\mathbf{S}_{\rho}$\}},\\
 \ker \nabla^{1,0}&=\textrm{\{anti-holomorphic sections of $\mathbf{S}_{\rho}$\}}.
 \end{split}\label{eqn:8-1}
\end{equation}
%%%%%%%%%%%% proposition %%%%%%%%%%%%%%%%
\begin{prop}\label{prop:8-1}
\begin{enumerate}
	\item The holomorphic (anti-holomorphic) sections of $\mathbf{S}_{\rho}$ on a closed K\"ahler manifold is characterized as an intersection of kernels of the K\"ahlerian gradients. 
\begin{equation}
\begin{split}
H^0(M,\mathbf{S}_{\rho})=\textrm{\{holomorphic sections of $\mathbf{S}_{\rho}$\}} &=\bigcap_{1\le i \le m} \ker D_{+i},\\
\textrm{\{anti-holomorphic sections of $\mathbf{S}_{\rho}$\}}&=\bigcap_{1\le i\le m}\ker D_{-i}.
\end{split}\label{eqn:8-2}
\end{equation}
Furthermore, a holomorphic section $\phi$ of $\mathbf{S}_{\rho}$ satisfies 
\begin{equation}
R^q_{\rho} (\phi) =\sum_i w_{-i}^q D_{-i}^{\ast}D_{-i}\phi, \label{holo-1}
\end{equation}
and an anti-holomorphic section $\phi$ does 
\begin{equation}
\tilde{R}^q_{\rho}(\phi)=\sum_i w_{+i}^q D_{+i}^{\ast}D_{+i}\phi.
\end{equation}
\item On a closed Ricci-flat K\"ahler manifold, 
\begin{equation}
\begin{split}
 &\textrm{\{holomorphic sections of $\mathbf{S}_{\rho}$\}}\\
=&\textrm{\{anti-holomorphic sections of $\mathbf{S}_{\rho}$\}}\\
=&\textrm{\{parallel sections of $\mathbf{S}_{\rho}$\}}.
\end{split}\label{eqn:8-4}
\end{equation}
\item If the Ricci curvature is non-positive and negative at some point, and $\pi_{\rho}(c_1)=\sum \rho^i$ is positive, then $H^0(M,\mathbf{S}_{\rho})=\{0\}$. On the other hand, if the Ricci curvature is non-negative and positive at some point, and $\pi_{\rho}(c_1)=\sum \rho^i$ is negative, then $H^0(M,\mathbf{S}_{\rho})=\{0\}$.
\end{enumerate}
\end{prop}
%%%%%%%%%
%%%%%%%%%%%%%   proof %%%%%%%%%%%%%%%%%%%
\begin{proof}
The first and second assertions follow from the Bochner identities in theorem \ref{thm:7-4}. If we assume that the Ricci curvature is non-positive and negative at some point, and $\pi_{\rho}(c_1)=\sum \rho^i$ is positive, then we show from \eqref{eqn:7-13} that $R^0_{\rho}$ is non-positive and negative at some point. For a holomorphic section $\phi$ of $\mathbf{S}_{\rho}$, we have 
$$
0 \le \sum \| D_{-i} \phi \|^2 = \int_M (R^0_{\rho} \phi,\phi) dx \le 0
$$
So we know that $\phi$ is parallel. Since $R^0_{\rho}$ is negative at some point, $\phi$ is zero section. 
\end{proof}
%%%%%%%
%%%%%%%%%%%%  example %%%%%%%%%%%%%%%
\begin{exam}[holomorphic vector fields]\label{ex:8-1}
The holomorphic vector fields are holomorphic sections of $T^{1,0}(M)=\mathbf{S}_{\mu_1}$, where $\mu_1$ is $(1,0_{m-1})$. So, if the Ricci curvature is non-positive and negative at some point, then $H^0(M,T^{1,0}(M))=\{0\}$. 
\end{exam}
%%%%%
%%%%%%%%%% example %%%%%%%%%%%%%%%
\begin{exam}[holomorphic differential forms]\label{ex:8-2}
The holomorphic $q$-forms are holomorphic sections of $\Lambda^{q,0}(M)=\mathbf{S}_{\rho}$, where $\rho$ is $(0_{m-q},(-1)_q)$. So, if the Ricci curvature is non-negative and positive at some point, then $H^0(M,\Lambda^{q,0}(M))=H^0(M,\Omega^q)=\mathbf{H}^{q,0}=\{0\}$ for $q\ge 1$.

\end{exam}
%%%%%%%
%%%%%%%%   example %%%%%%%%%%%%%%%%
\begin{exam}[holomorphic sections on $\mathbb{C}P^m$]\label{ex:8-3}
Let $M$ be a closed K\"ahler manifold of constant holomorphic sectional curvature $=r>0$, which can be identified with the complex projective space $\mathbb{C}P^m$ (see \cite{KN}). We take a non-zero holomorphic section $\phi$ of $\mathbf{S}_{\rho}$. Then it is from \eqref{holo-1} and \eqref{eqn:7-14} that 
$$
\frac{r}{2}\pi_{\rho}(c_qc_1+c_{q+1})\phi=\sum w_{-i}^q D_{-i}^{\ast}D_{-i}\phi.$$
Then we have 
\begin{equation}
D_{-i}^{\ast}D_{-i}\phi=\frac{r}{2}\gamma_{-i}(w_{-i}+\pi_{\rho}(c_1))\phi=\frac{r}{2}\gamma_{-i}(w_{-i}+\sum \rho^i)\phi. \nonumber
\end{equation}
Thus non-zero holomorphic sections of $\mathbf{S}_{\rho}$ on $\mathbb{C}P^m$ are eigensections for $D_{-i}^{\ast}D_{-i}$ with eigenvalue $r\gamma_{-i}(w_{-i}+\sum \rho^i)/2$ .
\end{exam}
%%%%%%%%

%%%%%%% subsection %%%%%%%%%%%%%%%
\subsection{The Bochner-Weitzenb\"ock formula}
We discuss the (generalized) Bochner-Weitzenb\"ock formula for the K\"ahlerian gradients. We know that there always exist $\mathbf{S}_{\rho+\mu_1}$ and $\mathbf{S}_{\rho-\mu_m}$ in the decomposition $\mathbf{S}_{\rho}\otimes T^{1,0}(M)=\sum \mathbf{S}_{\rho+\mu_i}$ and $\mathbf{S}_{\rho}\otimes T^{0,1}(M)=\sum \mathbf{S}_{\rho-\mu_i}$. So we always have the operators 
\begin{gather}
D_{+1}:\Gamma (\mathbf{S}_{\rho})\to \Gamma(\mathbf{S}_{\rho+\mu_1}),\label{eqn:8-5} \\ 
D_{-m}:\Gamma (\mathbf{S}_{\rho})\to \Gamma(\mathbf{S}_{\rho-\mu_m}).\label{eqn:8-6} 
\end{gather}
We call $D_{+1}$ \textit{the top operator} and $D_{-m}$ \textit{the bottom operator}. If we cancel the top and bottom operators from the Bochner identities with degree $0$ and $1$, we have the Bochner-Weitzenb\"ock formula and a vanishing theorem for the K\"ahlerian gradients.
%%%%%%%%%% proposition %%%%%%%%%%%%%%%
\begin{prop}
Let $\mathbf{S}_{\rho}$ be the associated bundle with rank $\ge 2$. Then the K\"ahlerian gradients $D_{\pm i}$ and $D^{\ast}_{\pm i}$ associated to $\rho$ satisfy that 
\begin{multline}
\sum_{1\le i\le m-1}\frac{2(w_{-m}-w_{-i})}{w_{+1}+w_{-m}}D_{-i}^{\ast}D_{-i}+\sum_{2\le i\le m}\frac{2(w_{+1}-w_{+i})}{w_{+1}+w_{-m}}D_{+i}^{\ast}D_{+i}\\=\nabla^{\ast}\nabla-\frac{2}{w_{+1}+w_{-m}}R_{\rho}^1-\frac{w_{+1}-w_{-m}}{w_{+1}+w_{-m}}R^0_{\rho}.\label{eqn:8-7} 
\end{multline}
If the curvature endomorphism 
$$-\frac{2}{w_{+1}+w_{-m}}R_{\rho}^1-\frac{w_{+1}-w_{-m}}{w_{+1}+w_{-m}}R^0_{\rho}$$
 is non-negative and positive at some point, then 
\begin{equation}
\bigcap_{1\le i\le m-1} \ker D_{-i}\cap\bigcap_{2\le i\le m}\ker D_{+i}=\{0\}.\nonumber
\end{equation}
\end{prop}
%%%%%%%%%%
We rewrite \eqref{eqn:8-7} with $\rho=(\rho^1,\cdots, \rho^m)$ and have
\begin{multline}
\sum_{1\le i\le m-1}\frac{2(\rho^i-\rho^m+m-i)}{\rho^1-\rho^m}D_{-i}^{\ast}D_{-i}+\sum_{2\le i\le m}\frac{2(\rho^1-\rho^i+i-1)}{\rho^1-\rho^m}D_{+i}^{\ast}D_{+i}\\=\nabla^{\ast}\nabla+\frac{2}{\rho^1-\rho^m}R_{\rho}^1-\frac{\rho^1+\rho^m}{\rho^1-\rho^m}R^0_{\rho}.\label{eqn:8-8}
\end{multline}
Since $\rho^1\ge \cdots\ge \rho^m$, the left hand side on \eqref{eqn:8-7} gives a non-negative operator. When the rank of $\mathbf{S}_{\rho}$ is $1$, the representation $(\pi_{\rho},V_{\rho})$ is a one-dimensional representation of $U(m)$ whose highest weight is $\rho=(k_m)$ for a $k$ in $\mathbb{Z}$. In this case, the above corollary is not valid because of $\rho^1=\rho^m$. So we have to assume that the rank of $\mathbf{S}_{\rho}$ is more than $1$ in the above proposition.

We discuss the case that the highest weight $\rho$ is $(k_m)$ for $k\in \mathbb{Z}$. We denote the canonical line bundle $\Lambda^{m,0}(M)$ of $M$ by  $K$, which is the associated bundle $\mathbf{S}_{((-1)_m)}$. So the line bundle $\mathbf{S}_{(k_m)}$ is isomorphic to $K^{-k}:=(K^{\ast})^k=(\bar{K})^k$. On this line bundle, we have two operators 
$$
D_{-m}:\Gamma(K^{-k})\to \Gamma (K^{-k}\otimes \Lambda^{1,0}(M)), \quad D_{+1}=\bar{\partial}:\Gamma(K^{-k})\to \Gamma (K^{-k}\otimes \Lambda^{0,1}(M)).
$$
There are two identities for these operators,
$$
D_{-m}^{\ast}D_{-m}+D_{+1}^{\ast}D_{+1}=\nabla^{\ast}\nabla, \quad D_{-m}^{\ast}D_{-m}-D_{+1}^{\ast}D_{+1}=\frac{k}{2}\kappa.
$$
Because of $D_{+1}=\bar{\partial}$, it holds that $2\bar{\partial}^{\ast} \bar{\partial}=\nabla^{\ast}\nabla-k\kappa/2$. 

%%%%%%%% example %%%%%%%%%%%%%%%%%%%%%%%%%%%%%%%%%%%%%%%%%%%%%%%%%%%%%%%%%%%%%%%%%%%%%%%%%%%%%%%%%%%%%%%%%%%%%%
\subsection{The Dolbeault-Dirac operator}
On a K\"ahler manifold $M$, there is a natural spin-$c$ structure and the spinor bundle $\sum \Lambda^{0,p}(M)$. From local calculations for spinors in section \ref{sec:5}, we can have some properties of the Dolbeault-Dirac operator on $\sum \Lambda^{0,p}(M)$.

On each vector bundle $\mathbf{S}_{\Lambda^{0,p}}=\Lambda^{0,p}(M)$, we have four operators,
\begin{equation}
\begin{split}
D_{-m}:& \Gamma (\Lambda^{0,p}(M))\to \Gamma (\mathbf{S}_{\Lambda^{0,p}-\mu_m}),\\ 
D_{-p}:& \Gamma (\Lambda^{0,p}(M))\to \Gamma (\mathbf{S}_{\Lambda^{0,p}-\mu_p})=\Gamma(\Lambda^{0,p-1}(M)),\\ 
D_{+1}:& \Gamma (\Lambda^{0,p}(M))\to 
         \Gamma (\mathbf{S}_{\Lambda^{0,p}+\mu_1}),\\
D_{+(p+1)}:& \Gamma (\Lambda^{0,p}(M))\to \Gamma (\mathbf{S}_{\Lambda^{0,p}+\mu_{p+1}})=\Gamma (\Lambda^{0,p+1}(M)).
\end{split}\label{eqn:8-9}
\end{equation}
Here, we know that 
\begin{gather}
\sqrt{m-p+1}D_{-p}=-\sum \bar{\e}_k\cdot \nabla_{\e_k}=-\bar{\partial}^{\ast}, \label{eqn:8-12}\\
\sqrt{p+1}D_{+(p+1)}=\sum \e_k\nabla_{\bar{\e}_k}=\bar{\partial}. \label{eqn:8-13}
\end{gather}
The first order differential operator $\sqrt{2}(\bar{\partial}+\bar{\partial}^{\ast})$ on $\sum \Lambda^{0,p}(M)$ is called \textit{the Dolbeault-Dirac operator}. 

It is from the Bochner identities and \eqref{eqn:cur} that
\begin{equation}
\begin{split}
D_{-m}^{\ast}D_{-m}+D_{-p}^{\ast}D_{-p}+D_{+1}^{\ast}D_{+1}+D_{+(p+1)}^{\ast}D_{+(p+1)}=\nabla^{\ast}\nabla,  \\
D_{-m}^{\ast}D_{-m}+D_{-p}^{\ast}D_{-p}-D_{+1}^{\ast}D_{+1}-D_{+(p+1)}^{\ast}D_{+(p+1)}=R_{\Lambda^{0,p}}^0, \\
(m-p+1)D_{-p}^{\ast}D_{-p}-D_{+1}^{\ast}D_{+1}+pD_{+(p+1)}^{\ast}D_{+(p+1)}=R_{\Lambda^{0,p}}^1=R_{\Lambda^{0,p}}^0. 
\end{split}\label{eqn:8-10}
\end{equation}
Then we have the Bochner-Weitzenb\"ock formula 
\begin{equation}
\begin{split}
 &2(m-p+1)D_{-p}^{\ast}D_{-p}+2(p+1)D_{+(p+1)}^{\ast}D_{+(p+1)} \\
=&\nabla^{\ast}\nabla+R_{\Lambda^{0,p}}^0  =2(\bar{\partial}\bar{\partial}^{\ast}+\bar{\partial}^{\ast}\bar{\partial}).
\end{split}\label{eqn:8-11}
\end{equation}
On the other hand, because of $\nabla^{1,0}=\partial$ on $\Lambda^{0,p}(M)$, we have 
\begin{equation}
D_{-m}+D_{-p}=\partial. \label{eqn:8-14}
\end{equation}
The following proposition is well-known in K\"ahler geometry. 
%%%%%%%% proposition %%%%%%%%%%%%%%
\begin{prop}
\begin{enumerate}
	\item If the Ricci curvature of $M$ is non-negative and positive at some point, then 
\begin{equation}
\dim \mathbf{H}^{0,p}=\dim \ker \bar{\partial}|_{\Lambda^{0,p}} \cap \ker \bar{\partial}^{\ast}|_{\Lambda^{0,p}}=0 \label{eqn:8-15}
\end{equation}
for $p=1,\cdots,m$. 
\item A harmonic $(0,p)$-form $\phi$ satisfies 
\begin{equation}
\partial\phi=0, \quad \nabla^{\ast}\nabla \phi=R_{\Lambda^{0,p}}\phi.
\end{equation}
\item A holomorphic $(0,p)$-form $\phi$ satisfies 
\begin{equation}
\partial^{\ast}\partial\phi=\bar{\partial}\bar{\partial}^{\ast}\phi=-R_{\Lambda^{0,p}}^0\phi, \label{eqn:8-15-1}
\end{equation}
where we say that $\phi$ is holomorphic if $\nabla^{0,1}\phi=0$.
\item If $\mathrm{Ric}\ge c>0$, then every eigenvalue of $2(\bar{\partial}\bar{\partial}^{\ast}+\bar{\partial}^{\ast}\bar{\partial})$ on $\Lambda^{0,p}(M)$ satisfies\begin{equation}
\lambda \ge 2pc. \label{eqn:8-15-2}
\end{equation} 
\end{enumerate}
\end{prop}
%%%%%%%%
%%%%%%%%%%% proof %%%%%%%%%%
\begin{proof}
If the Ricci curvature $\mathrm{Ric}$ is non-negative and positive at some point, then $R_{\Lambda^{0.p}}^0$ is non-negative and positive at some point. The equation \eqref{eqn:8-15} follows from the Bochner-Weitzenb\"ock formula \eqref{eqn:8-11}. The second and third claims follow from \eqref{eqn:8-10}. We shall prove the fourth claim. It follows from \eqref{eqn:8-10} that, for a $(0,p)$-eigenform $\phi$ with eigenvalue $\lambda$, 
\begin{equation}
\begin{split}
 \lambda \|\phi\|^2=&2(m-p+1)\|D_{-p}\phi\|^2+2(p+1)\|D_{+(p+1)}\phi\|^2 \\
\ge & 2(m-p+1)\|D_{-p}\phi\|^2+2p\|D_{+(p+1)}\phi\|^2 \\
=&2\|D_{+1}\phi\|^2+2\int_M (R_{\Lambda^{0,p}}^0\phi,\phi) \\
 \ge & 2\int_M (R_{\Lambda^{0,p}}^0\phi,\phi) \ge 2cp\|\phi\|^2.
\end{split}\nonumber
\end{equation}
We know that, if the equality sign holds on \eqref{eqn:8-15-2}, then $\phi$ is holomorphic. 
\end{proof}
%%%%%%%%%
%%%%%%%%%%%%  example %%%%%%%%%%%%%%%%%
\subsection{The Dirac operator}
We discuss the Dirac operator on a spin K\"ahler manifold. A K\"ahler manifold $M$ has a spin structure if and only if there exists a holomorphic square root $L:=\sqrt{K}$ of the canonical line bundle $K=\Lambda^{m,0}(M)$. The spinor bundle on a spin K\"ahler manifold is $\mathbf{S}=\sum \Lambda^{0,p}(M)\otimes L=\sum \Lambda^{0,p}(M)\otimes \sqrt{K}$. It is well-known that there is an anti-linear bundle isomorphism 
\begin{equation}
\mathbf{j}:\Lambda^{0,p}(M)\otimes L\ni \psi \mapsto \mathbf{j}(\psi) \in \Lambda^{0,m-p}(M)\otimes L=\overline{\Lambda^{0,p}(M)\otimes L}
\end{equation}
such that $\mathbf{j}$ is parallel, and commute or anti-commute with Clifford multiplications of real vectors, and $\mathbf{j}^2=(-1)^{m(m+1)/2}$ (cf. \cite{K4}). 

On this spinor bundle, we have four operators, 
\begin{equation}
\begin{split}
D_{-m}^L:& \Gamma (\Lambda^{0,p}(M)\otimes L)\to \Gamma (\mathbf{S}_{\Lambda^{0,p}-\mu_m}\otimes L),\\ 
D_{-p}^L:& \Gamma (\Lambda^{0,p}(M)\otimes L)\to \Gamma(\Lambda^{0,p-1}(M)\otimes L), \\ 
D_{+1}^L:& \Gamma (\Lambda^{0,p}(M)\otimes L)\to 
         \Gamma (\mathbf{S}_{\Lambda^{0,p}+\mu_1}\otimes L),\\
D_{+(p+1)}^L:& \Gamma (\Lambda^{0,p}(M)\otimes L)\to \Gamma (\Lambda^{0,p+1}(M)\otimes L).
\end{split}
\end{equation}
Here, we remark that, with respect to $\mathbf{j}$, the operators $D_{-m}^L$, $D_{+1}^L$, $D_{-p}^L$, and $D_{+(p+1)}^L$ on $\Lambda^{0,p}(M)\otimes L$ correspond to  $D_{+1}^L$, $D_{-m}^L$, $D_{+(m-p+1)}^L$, and $D_{-(m-p)}^L$ on $\Lambda^{0,m-p}(M)\otimes L$, respectively.

Now, we calculate the curvature endomorphisms $\mathfrak{R}_L^0$ and $\mathfrak{R}^1_L$. The curvature of $L$ is 
$$
R_L(\e_k,\bar{\e}_l)=R_{\sqrt{\Lambda^{m,0}}}(\e_k,\bar{\e}_l)=\frac{1}{2} R_{\Lambda^{m,0}}(\e_k,\bar{\e}_l)=-\frac{1}{2}\sum g(R_T(\e_i,\bar{\e}_i)\e_k,\bar{\e}_l).
$$
Then we have 
\begin{gather}
\mathfrak{R}^0_L=\id \otimes \sum R_{L}(\e_k,\bar{\e}_k)=\id \otimes \frac{1}{2}R_{\Lambda^{m,0}}^0=-\frac{1}{4}\kappa, \nonumber \\
\mathfrak{R}^1_L=\sum \pi_{\Lambda^{0,p}}(e_{kl})\otimes R_{L}(\e_l,\bar{\e}_k)=-\frac{1}{2}R_{\Lambda^{0,p}}^1=-\frac{1}{2}R_{\Lambda^{0,p}}^0. \nonumber
\end{gather}
It is from the Bochner identities that
%%%%%%  Bochner identities for Dirac %%%%%%%%%%
\begin{equation}
\begin{split}
(D_{-m}^L)^{\ast}D_{-m}^L+(D_{-p}^L)^{\ast}D_{-p}^L+(D_{+1}^L)^{\ast}D_{+1}^L+(D^L)_{+(p+1)}^{\ast}D^L_{+(p+1)}=\nabla^{\ast}\nabla,  \\
(D_{-m}^L)^{\ast}D_{-m}^L+(D_{-p}^L)^{\ast}D_{-p}^L-(D_{+1}^L)^{\ast}D^L_{+1}-(D^L)_{+(p+1)}^{\ast}D^L_{+(p+1)} \\
 =R_{\Lambda^{0,p}}^0+\mathfrak{R}_L^0=R_{\Lambda^{0,p}}^0-\frac{1}{4}\kappa, \\  (m-p+1)(D_{-p}^L)^{\ast}D_{-p}^L-(D_{+1}^L)^{\ast}D_{+1}^L+p(D_{+(p+1)}^L)^{\ast}D_{+(p+1)}^L \\
 =R_{\Lambda^{0,p}}^1+\mathfrak{R}_L^1 =\frac{1}{2}R^0_{\Lambda^{0,p}}. 
\end{split}\label{eqn:8-17}
\end{equation}
%%%%%%%%%%%%
The Bochner-Weitzenb\"ock formula gives the Lichnerowicz formula for the Dirac operator $D:=\sqrt{2}(\bar{\partial}_L+\bar{\partial}_L^{\ast})$,
%%%%%%%%%%%% Lichnerowicz formula for the Dirac  %%%%%%%%%
\begin{equation}
\begin{split}
 &2(m-p+1)(D_{-p}^L)^{\ast}D_{-p}^L+2(p+1)(D_{+(p+1)}^L)^{\ast}D_{+(p+1)}^L \\
=&\nabla^{\ast}\nabla+\frac{1}{4}\kappa   =2(\bar{\partial}_L\bar{\partial}_L^{\ast}+\bar{\partial}_L^{\ast}\bar{\partial}_L)=D^2
\end{split}\label{eqn:8-18}
\end{equation}
%%%%%%%
on $\Lambda^{0,p}(M)\otimes L$ for any $p$. 

The Bochner identities are applied to the following proposition shown by K. D. Kirchberg in \cite{K4}.
%%%%%%%%%%%%%%%%%%%%%%%%  prop %%%%%%%%%%%%%%%%%%%%%%%
\begin{prop}[\cite{K4}]
\begin{enumerate}
\item A holomorphic spinor $\phi$ of $\Lambda^{0,p}(M)\otimes L$ satisfies that 
\begin{equation}
\bar{\partial}_L\bar{\partial}_L^{\ast}\phi=\frac{1}{2}R^0_{\Lambda^{0,p}}\phi.\label{eqn:8-19}
\end{equation}
\item If $R_{\Lambda^{0,p}}^0-\kappa/4$ is non-positive and negative at some point, there is no holomorphic spinor of $\Lambda^{0,p}(M)\otimes L$. 
\item Let $M$ be a closed spin K\"ahler-Einstein manifold. If $\kappa>0$ and $p\le m/2$, then there is no holomorphic spinor of $\Lambda^{0,p}(M)\otimes L$.
\end{enumerate}
\end{prop}
%%%%%%%%
%%%%%%%%%%%%% proof %%%%%%%%%%%%%%%%%%
\begin{proof}
The first and second assertions easily follow from \eqref{eqn:8-17}. We shall prove the third one. We have known that $R^0_{\Lambda^{0,m}}=\frac{1}{2}\kappa$, where $\kappa$ is constant on a spin K\"ahler-Einstein manifold. So we have
$$
R_{\Lambda^{0,p}}^0-\frac{1}{4}\kappa=R^0_{\Lambda^{0,p}}-\frac{1}{2}R^0_{\Lambda^{0,m}}=(\frac{p}{2m}-\frac{1}{4})\kappa.
$$
If $\kappa>0$ and $p< m/2$, then $R_{\Lambda^{0,p}}^0-\kappa/4$ is negative and there is no holomorphic spinor. In the case that $\kappa>0$ and $p=m/2$, we know that a holomorphic spinor is parallel spinor. Since the Ricci curvature is zero on a spin manifold with parallel spinor, we have a contradiction to $\kappa>0$. 
\end{proof}
%%%%%%%%%%%%

We investigate the operators $D^L_{-m}$ and $D^L_{+1}$. These operators are known as the K\"ahlerian twistor operators on $M$. It follows from \eqref{eqn:5-6} that 
\begin{equation}
\begin{split}
 &\nabla^{0,1}\phi\\
=&\sum \nabla_{\bar{\e}_k}\phi\otimes \e_k \\
=&-\frac{1}{p+1}\sum \bar{\e}_i\cdot\e_k\cdot\nabla_{\bar{\e}_k}\phi\otimes \e_i+(\sum \nabla_{\bar{\e}_k}\phi\otimes \e_k+\frac{1}{p+1}\sum \bar{\e}_i\cdot\e_k\cdot\nabla_{\bar{\e}_k}\phi\otimes \e_i)\\
=&-\frac{1}{p+1}\sum \bar{\e}_i\cdot\bar{\partial}_L\phi\otimes \e_i+\sum (\nabla_{\bar{\e}_i}\phi+\frac{1}{p+1}\sum \bar{\e}_i\cdot\bar{\partial}_L\phi)\otimes \e_i.
\end{split}\nonumber
\end{equation}
Similarly, we show from \eqref{eqn:5-7} that
\begin{equation}
\nabla^{1,0}\phi=-\frac{1}{m-p+1}\sum \e_i\cdot\bar{\partial}_L^{\ast}\phi\otimes \bar{\e}_i+\sum (\nabla_{\e_i}\phi+\frac{1}{m-p+1}\sum \e_i\cdot\bar{\partial}_L^{\ast}\phi)\otimes \bar{\e}_i. \nonumber
\end{equation}
%%%%%%%%%%  prop %%%%%%%%%%%%%%%%%%%%%%%%
\begin{prop}
Let $\phi$ be a spinor in $\Gamma (\Lambda^{0,p}(M)\otimes L)$. Then, 
\begin{gather}
\phi\in\ker D_{-m}^L \iff \nabla_{\bar{X}}\phi+\frac{1}{p+1}\bar{X}\cdot\bar{\partial}_L\phi=0 \quad \textrm{for any $X$ in $\Gamma (T^{1,0}(M))$},\label{eqn:8-20}\\
\phi\in\ker D_{+1}^L \iff \nabla_{X}\phi+\frac{1}{m-p+1} X\cdot\bar{\partial}^{\ast}_L\phi=0 \quad \textrm{for any $X$ in $\Gamma (T^{1,0}(M))$}.\label{eqn:8-21}
\end{gather}
The spinor $\phi$ is called a (K\"ahlerian) twistor spinor if $\phi$ is in $\ker D_{-m}^L\cap \ker D_{+1}^L$.
\end{prop}
%%%%%%%%%%%%%

In \cite{K1} and \cite{K2}, K. D. Kirchberg has given an estimate for the eigenvalues of $D^2$ on a spin K\"ahler manifold. We shall give an another proof of the estimate by using our Bochner identities. It is from the Bochner identities \eqref{eqn:8-17} on the spinor bundle that 
\begin{equation}
\begin{cases}
2\bar{\partial}_L^{\ast}\bar{\partial}_L=(D_{-m}^L)^{\ast}D_{-m}^L+\frac{1}{4}\kappa & \textrm{for $p=0$}, \\
2\bar{\partial}_L\bar{\partial}_L^{\ast}=(D_{+1}^L)^{\ast}D_{+1}^L+\frac{1}{4}\kappa & \textrm{for $p=m$}, \\
 \frac{2m-2p+1}{2m-2p+2}2\bar{\partial}_L\bar{\partial}^{\ast}_L+\frac{2p+1}{2p+2}2\bar{\partial}^{\ast}_L\bar{\partial}_L =(D_{-m}^L)^{\ast}D_{-m}^L+(D_{+1}^L)^{\ast}D_{+1}^L+\frac{1}{4}\kappa  &\textrm{otherwise}.
\end{cases}\label{eqn:8-22}
\end{equation}
If we use only these equations, then we have the Friedrich's eigenvalue estimate of the Dirac operator, that is, $\lambda^2\ge \frac{1}{4}\frac{2m}{2m-1}\min_x \kappa(x)$. To more sharp estimate on a spin K\"ahler manifold,  we need the Hodge-de Rham-Kodaira decomposition for $L$,
\begin{equation}
\Gamma (\Lambda^{0,p}(M)\otimes L)=\mathbf{H}^{0,p}(L)\oplus \bar{\partial}_L\Gamma (\Lambda^{0,p}(M)\otimes L)\oplus \bar{\partial}_L^{\ast}\Gamma (\Lambda^{0,p}(M)\otimes L).\label{eqn:8-23}
\end{equation}
This decomposition gives the following well-known fact.
%%%%%%%%%%%% lemma %%%%%%%%%%%%%%%%%%%%
\begin{lem}
Let $\lambda\neq 0$ be an eigenvalue of the Dirac operator $D=\sqrt{2}(\bar{\partial}_L+\bar{\partial}_L^{\ast})$ and $E_{\lambda}$ be the eigenspace with eigenvalue $\lambda$. Then there is a decomposition $E_{\lambda}=\sum_{p=0}^{m-1} E_{\lambda}^p$ such that any spinor $\Psi_p$ in $E_{\lambda}^p$ is given by 
\begin{equation}
\begin{split}
\Psi_p=\psi_p+\psi_{p+1}\in \Gamma (\Lambda^{0,p}(M)\otimes L)\oplus \Gamma (\Lambda^{0,p+1}(M)\otimes L), \\
\sqrt{2}\bar{\partial}_L\psi_p=\lambda\psi_{p+1}, \quad \sqrt{2}\bar{\partial}^{\ast}_L\psi_{p+1}=\lambda\psi_{p}, \quad \bar{\partial}^{\ast}_L\psi_p=0, \quad \bar{\partial}_L\psi_{p+1}=0.
\end{split}\label{eqn:8-24}
\end{equation}
\end{lem}
%%%%%%
%%%%%%%%%%%%%% proof %%%%%%%%%%%%%%%%%%%%%%%
\begin{proof}
Let $\eta$ be an eigenspinor with eigenvalue $\lambda\neq 0$ for the Dirac operator $D$. We decompose $\eta$ as 
$$
\eta=\sum_{p=0}^m \eta_p, \quad \eta_p \in \Gamma (\Lambda^{0,p}(M)\otimes L).
$$
Since $D(\lambda^{-2}D\eta_p)=\eta_p$, we have
\begin{equation}
\eta_p=\bar{\partial}_L \phi_{p,p-1}+\bar{\partial}^{\ast}_L\phi_{p,p+1}, \nonumber
\end{equation}
where we set 
$$
\phi_{p,p-1}:=2\lambda^{-2}\bar{\partial}_L^{\ast}\eta_p \in \Gamma(\Lambda^{0,p-1}(M)\otimes L), \quad \phi_{p,p+1}:=2\lambda^{-2}\bar{\partial}_L\eta_p\in\Gamma(\Lambda^{0,p+1}(M)\otimes L)
$$ and remark that $\phi_{0,-1}=\phi_{m,m+1}=0$. Since $\eta=\sum \eta_p$ is the eigenspinor with eigenvalue $\lambda$, it holds that 
\begin{equation}
\sqrt{2}\bar{\partial}_L\bar{\partial}_L^{\ast}\phi_{p-1,p}+\sqrt{2}\bar{\partial}_L^{\ast}\bar{\partial}_L\phi_{p+1,p}=\lambda (\bar{\partial}_L\phi_{p,p-1}+\bar{\partial}_L^{\ast}\phi_{p,p+1})
\nonumber
\end{equation}
for $p=0,\cdots,m$. We use the Hodge-de Rham-Kodaira decomposition on this equation and have
\begin{equation}
 \sqrt{2}\bar{\partial}_L^{\ast}\bar{\partial}_L\phi_{p+1,p}=\lambda \bar{\partial}_L^{\ast} \phi_{p,p+1},   \quad \sqrt{2}\bar{\partial}_L\bar{\partial}_L^{\ast}\phi_{p,p+1}=\lambda \bar{\partial}_L\phi_{p+1,p} \nonumber
\end{equation}
for $p=0,\cdots,m$. Define
\begin{equation}
\Psi_p :=\bar{\partial}^{\ast}_L \phi_{p,p+1}+\bar{\partial}_L \phi_{p+1,p} \quad \textrm{ for $p=0,\cdots,m-1$} , \nonumber 
\end{equation}
and each $\Psi_p$ is an eigenspinor of the Dirac operator with eigenvalue $\lambda$ and $\eta=\sum_{p=0}^{m-1} \Psi_p$. 
\end{proof}
%%%%%%%%%%%%%%%
By this lemma,  we consider only the eigenspinor $\Psi_p=\psi_p+\psi_{p+1}$ satisfying the condition \eqref{eqn:8-24}. First, we apply the identity \eqref{eqn:8-22} to $\psi_p$ and have 
\begin{equation}
\begin{split}
\frac{2p+1}{2(p+1)}\lambda^2\|\psi_p\|^2 &=\|D_{-m}^L\psi_p\|^2+\|D_{+1}^L\psi_p\|^2+\frac{1}{4}\int_M (\kappa (x)\psi_p,\psi_p)dx \\
  &\ge \frac{1}{4}\int_M (\kappa (x)\psi_p,\psi_p)dx \ge \frac{\kappa_0}{4}\|\psi_p\|^2,
\end{split}\nonumber
\end{equation}
where $\kappa_0=\min_{x\in M}\kappa (x)$. So the eigenvalue $\lambda$ has a lower bound depending on the scalar curvature $\kappa$,
\begin{equation}
\lambda^2 \ge \frac{2p+2}{2p+1}\frac{\kappa_0}{4}. \nonumber
\end{equation}
On the other hand, we apply the same identity \eqref{eqn:8-22} to $\psi_{p+1}$ and have 
\begin{equation}
\begin{split}
  &\frac{2m-2(p+1)+1}{2(m-(p+1)+1)}\lambda^2\|\psi_{p+1}\|^2 \\
=&\|D_{-m}^L\psi_{p+1}\|^2+\|D_{+1}^L\psi_{p+1}\|^2+\frac{1}{4}\int_M (\kappa (x)\psi_{p+1},\psi_{p+1})dx 
  \ge \frac{\kappa_0}{4}\|\psi_{p+1}\|^2.
\end{split}\nonumber
\end{equation}
Then the eigenvalue $\lambda$ has another lower bound, 
\begin{equation}
\lambda^2\ge \frac{2m-2p}{2m-2p-1}\frac{\kappa_0}{4}. \nonumber
\end{equation}
Thus, the eigenspinor $\Psi_p$ with nonzero eigenvalue $\lambda$ has a lower bound
\begin{equation}
\lambda^2\ge \frac{\kappa_0}{4}\max \{\frac{2p+2}{2p+1},\frac{2m-2p}{2m-2p-1}\}.\end{equation}
Since this estimate holds for $p=0,\cdots,m-1$, we conclude that non-zero eigenvalue $\lambda$ for the Dirac operator satisfies 
\begin{equation}
\lambda^2\ge \frac{\kappa_0}{4}\min_{0\le p\le m-1}\max \{\frac{2p+2}{2p+1},\frac{2m-2p}{2m-2p-1}\}. \nonumber
\end{equation}

In the case that $m$ is even, we have 
$$
\lambda^2\ge \frac{\kappa_0}{4}\frac{m}{m-1}.
$$
Furthermore, we know that, if the equality sign holds on this inequality, then there exists a spinor $\psi_{m/2-1}$ in $\Gamma (\Lambda^{0,m/2-1}(M)\otimes L)$ such that 
\begin{equation}
\psi_{m/2-1} \in \ker D_{-m}^L\cap \ker D_{+1}^L\cap \ker \bar{\partial}^{\ast}_L, \label{twistor 1}
\end{equation}
or a spinor $\phi_{m/2+1}$ in $\Gamma (\Lambda^{0,m/2+1}(M)\otimes L)$ such that\begin{equation}
\psi_{m/2+1} \in \ker D_{-m}^L\cap \ker D_{+1}^L\cap \ker \bar{\partial}_L.
\label{twistor 2}
\end{equation}
We show that, if a spinor $\psi_{m/2-1}$ satisfying \eqref{twistor 1}, then $\mathbf{j}(\psi_{m/2-1})$ is a spinor satisfying \eqref{twistor 2}.

In the case that $m$ is odd, we have 
$$
\lambda^2\ge \frac{\kappa_0}{4}\frac{m+1}{m}.
$$
If the equality sign holds on this inequality, then there exists a spinor $\psi_{m/2-1/2}$ in $\Gamma (\Lambda^{0,m/2-1/2}(M)\otimes L)$ such that 
$$
\psi_{m/2-1/2} \in \ker D_{-m}^L\cap \ker D_{+1}^L\cap \ker \bar{\partial}^{\ast}_L,
$$
or a spinor $\psi_{m/2+1/2}$ in $\Gamma (\Lambda^{0,m/2+1/2}(M)\otimes L)$ such that 
$$
\psi_{m/2+1/2} \in \ker D_{-m}^L\cap \ker D_{+1}^L\cap \ker \bar{\partial}_L.
$$
%%%%%%%%%%%%%%% prop %%%%%%%%%%%%%%%%%%%%%%
\begin{prop}[\cite{K1}, \cite{K2}]
Every eigenvalue of the Dirac operator on a closed spin K\"ahler manifold with positive scalar curvature satisfies
\begin{equation}
\begin{cases}
 \lambda^2\ge \frac{\kappa_0}{4}\frac{m}{m-1} & \textrm{for $m=$even}, \\
 \lambda^2\ge \frac{\kappa_0}{4}\frac{m+1}{m} & \textrm{for $m=$odd},
\end{cases}
\end{equation}
where $\kappa_0=\min_x\kappa(x)$. If the equality sign holds on the above inequality, then the scalar curvature is constant and there exists a twistor spinor $\psi$ such that 
\begin{equation}
\begin{cases}
\psi \in \Gamma (\Lambda^{0,m/2+1}(M)\otimes L)   & \textrm{for $m=$even}, \\ 
\psi \in \Gamma (\Lambda^{0,m/2+1/2}(M)\otimes L)  & \textrm{for $m=$odd},
\end{cases} 
\end{equation}
and
\begin{equation}
\psi \in \ker D_{-m}^L\cap \ker D_{+1}^L \cap \ker \bar{\partial}_L.
\end{equation}
\end{prop}
%%%%%%%%%%%%%%%%%%%%%

%%%%%%%%%%%%%%%%%%%%%%%%%%%%%%%%
%%  Acknowledgements          %%
%%%%%%%%%%%%%%%%%%%%%%%%%%%%%%%%
\section*{Acknowledgements}\label{sec:acknowledgement}
The author is partially supported by the Grant-in-Aid for Scientific Research (No. 13740120) from the Ministry of Education, Culture, Sports, Science and Technology. He also thanks Jun-ichi Okuda and Tatsuo Suzuki for discussions. 
%%%%%%%%%%%%%%%%%%%%%%%%%%%%%%%%
%        reference             %
%%%%%%%%%%%%%%%%%%%%%%%%%%%%%%%%
%%

%%%%%%%%%%%%%%%%%%%%%%%%%
\end{document}